\documentclass[psamsfonts]{amsart}
\usepackage{amssymb,eucal,latexsym,graphics}

\theoremstyle{plain}

\newtheorem{lemma}{Lemma}[section]
\newtheorem{theorem}[lemma]{Theorem}
\newtheorem{proposition}[lemma]{Proposition}
\newtheorem{corollary}[lemma]{Corollary}

\newtheorem*{stat}{\name}
\newcommand{\name}{testing}

\theoremstyle{definition}
\newtheorem{definition}[lemma]{Definition}

\newtheorem{problem}{Problem}

\theoremstyle{remark}

\newenvironment{all}[1]{\renewcommand{\name}{#1}\begin{stat}}
                         {\end{stat}}

\newcommand{\qedc}{{\qed}~{\rm Claim~{\theclaim}.}}

\numberwithin{equation}{section}

\newcommand{\tvi}{\vrule height 12pt depth 6pt width 0pt}

\newcommand{\pup}[1]{\textup{(}{#1}\textup{)}}

\newcommand{\xa}{{\boldsymbol{a}}}
\newcommand{\xb}{{\boldsymbol{b}}}
\newcommand{\xc}{{\boldsymbol{c}}}

\newcommand{\xe}{{\boldsymbol{e}}}

\newcommand{\xu}{{\boldsymbol{u}}}
\newcommand{\xv}{{\boldsymbol{v}}}

\newcommand{\bp}{\mathbin{\square}}
\newcommand{\ltp}{\mathbin{\boxtimes}}
\newcommand{\bt}[1]{M_3\langle{#1}\rangle}

\newcommand{\URP}{\mathrm{URP}}
\newcommand{\WURP}{\mathrm{WURP}}
\newcommand{\CLP}{\textup{CLP}}

\newcommand{\ccb}{conditionally co-Brou\-wer\-ian}
\newcommand{\ckcb}{conditionally $\kappa$-co-Brou\-wer\-ian}

\newcommand{\eps}{\varepsilon}
\newcommand{\es}{\varnothing}

\newcommand{\onto}{\twoheadrightarrow}

\newcommand{\seq}[1]{\langle{#1}\rangle}
\newcommand{\famm}[2]{\left\langle#1\mid#2\right\rangle}

\newcommand{\set}[1]{\{#1\}}
\newcommand{\setm}[2]{\set{#1\mid#2}}

\newcommand{\ol}[1]{\overline{#1}}

\newcommand{\dd}{\mathrm{d}}

\newcommand{\go}{\omega}

\DeclareMathOperator{\Var}{Var}

\newcommand{\CL}{\mathbf{L}}
\newcommand{\CD}{\mathbf{D}}
\newcommand{\CS}{\mathbf{S}}
\newcommand{\CSd}{\mathbf{S_d}}
\newcommand{\CSfd}{\mathbf{S_{fd}}}
\newcommand{\CSfb}{\mathbf{S_{fb}}}

\newcommand{\PL}{\mathbf{PL}}
\newcommand{\CC}{\mathbf{C}}
\newcommand{\MM}{\mathbf{M}}
\newcommand{\UU}{\mathbf{U}}
\newcommand{\VV}{\mathbf{V}}
\newcommand{\QQ}{\mathbb{Q}}

\newcommand{\DD}{\mathcal{D}}
\newcommand{\EE}{\mathcal{E}}
\newcommand{\LL}{\mathcal{L}}
\newcommand{\PP}{\mathcal{P}}

\newcommand{\MPL}[1]{$#1$-meas\-ured par\-tial lat\-tice}
\newcommand{\ML}[1]{$#1$-meas\-ured lat\-tice}

\DeclareMathOperator{\rng}{rng}
\DeclareMathOperator{\Con}{Con}
\DeclareMathOperator{\Aut}{Aut}

\DeclareMathOperator{\Conc}{Con_c}
\DeclareMathOperator{\Idc}{Id_c}
\DeclareMathOperator{\M}{M}

\newcommand{\id}{\mathrm{id}}
\newcommand{\jz}{$\langle\vee,0\rangle$}
\newcommand{\jzu}{$\langle\vee,0,1\rangle$}
\newcommand{\jzs}{\jz-semi\-lat\-tice}
\newcommand{\jzus}{\jzu-semi\-lat\-tice}
\newcommand{\jzh}{\jz-ho\-mo\-mor\-phism}
\newcommand{\js}{join-sem\-i\-lat\-tice}
\newcommand{\jh}{join-ho\-mo\-mor\-phism}

\newcommand{\fin}[1]{[#1]^{<\go}}
\newcommand{\fine}[1]{[#1]_*^{<\go}}

\DeclareMathOperator{\Fg}{F}
\DeclareMathOperator{\Fb}{FB}

\newcommand{\FL}{\Fg_{\mathbf{L}}}

\begin{document}

\dedicatory{Dedicated to Ralph McKenzie on his 60-th birthday}

\title[Congruence lattices of lattices]%
{A survey of recent results on congruence
lattices of lattices}

  \author[J.~T\r uma]{Ji\v r\'\i\ T\r uma}
  \address{Department of Algebra\\
           Faculty of Mathematics and Physics\\
           Sokolovsk\'a 83\\
           Charles University\\
           186 00 Praha 8\\
           Czech Republic}
  \email{tuma@karlin.mff.cuni.cz}

\author[F.~Wehrung]{Friedrich Wehrung}
\address{CNRS, UMR 6139\\
D\'epartement de Math\'ematiques\\
Universit\'e de Caen\\
14032 Caen Cedex\\
France}
\email{wehrung@math.unicaen.fr}

\urladdr{http://www.math.unicaen.fr/\~{}wehrung}

\date{\today}



\thanks{The first author was partially supported by GA
UK grant no.~162/1999 and by GA CR grant no. 201/99. The second author
was partially supported by the Fund of Mobility of the Charles
University (Prague), by FRVS grant no. 2125, by institutional grant
CEZ:J13/98:113200007, and by the Barrande program}

\begin{abstract}
We review recent results on congruence lattices of (infinite) lattices.
We discuss results obtained with box products, as well
as categorical, ring-theoretical, and topological results.
\end{abstract}

\maketitle
\tableofcontents

\section{Introduction}\label{S:Intro}

For a lattice $L$, the \emph{congruence lattice of $L$}, denoted here 
by $\Con L$,
is the lattice of all congruences of $L$ under inclusion. As the 
congruence lattice
of any algebraic system, the lattice $\Con L$ is \emph{algebraic}. The compact
elements of $\Con L$ are the \emph{finitely generated congruences}, 
that is, the
congruences of the form
  \[
  \bigvee_{i<n}\Theta_L(a_i,b_i),
  \]
where $n<\omega$, $a_i$, $b_i\in L$, for all $i<n$, and
$\Theta_L(a_i,b_i)$ (the \emph{principal congruence} generated by the
pair $\seq{a_i,b_i}$) denotes the least congruence of $L$ that identifies
$a_i$ and $b_i$. We denote by $\Conc L$, the \emph{congruence
semilattice of $L$}, the \jzs\ of all compact congruences of $L$.
A classical result by N. Funayama and T.~Nakayama \cite{FuNa42} states that the
lattice $\Con L$ is \emph{distributive}. Hence the \jzs\ $\Conc L$
is distributive, that is, for all $\xa$, $\xb$, $\xc\in\Conc L$,
if $\xc\leq\xa\vee\xb$, then there are elements $\xa'\leq\xa$ and $\xb'\leq\xb$
such that $\xc=\xa'\vee\xb'$. Most of the concepts we shall use in the present
paper are more conveniently expressed with $\Conc$ than with $\Con$.

Since the congruence lattice of any algebra is an algebraic lattice, it follows
that the congruence lattice of any lattice is an algebraic 
distributive lattice.
The question whether the converse of this result holds, that is, whether any
algebraic distributive lattice is isomorphic to $\Con L$, for some lattice $L$,
has been raised in the early forties by R.\,P. Dilworth, who solved the finite
case. We refer to this problem as the \emph{Congruence Lattice 
Problem}, \CLP\ in
short. The semilattice formulation of \CLP\ asks whether every 
distributive \jzs\
is isomorphic to $\Conc L$, for some lattice $L$.

Since the problem was raised, much progress has been done; we refer the
reader to G. Gr\"atzer and E.\,T. Schmidt \cite{GrScC} for a pre-1998
survey. Furthermore, it turns out that the topic of congruence lattices of
lattices can be divided into two parts: congruence lattices of
\emph{finite} lattices, and congruence lattices of \emph{infinite} 
lattices. These
topics are nearly disjoint (surprisingly?), although there are a few noteworthy
interactions between the two of them. We refer the reader to G. Gr\"atzer and
E.\,T. Schmidt \cite{FCSurv} for a survey of congruence lattices of 
finite lattices.

About the infinite case, the last few years have seen the emergence of many new
techniques and results about \CLP\ that the present paper intends to
survey. The main ideas can be separated into different groups.

\begin{itemize}
\item\textbf{Uniform refinement properties} (Section~\ref{S:URP}). Most
known partial negative solutions to \CLP\ are obtained \emph{via} certain
infinitary sentences of the theory of semilattices that hold in all
semilattices of the form $\Conc L$, for~$L$ in large classes of 
lattices, such as
the class of relatively complemented lattices. On the other hand, 
these formulas
do not hold in all distributive \jzs s. It also turns out that all 
the presently
known representation theorems yield semilattices with the strongest known
`uniform refinement property', which we denote here by $\URP^+$ (see
Propositions~\ref{P:RelCplURP+} and \ref{P:WOlim}).

\item\textbf{The $\bt{L}$ construction, tensor product, and box product}
(Section~\ref{S:BoolTr}). It was proposed as an open problem in G. 
Gr\"atzer and
E.\,T. Schmidt \cite{GrSc95}, whether every nontrivial lattice has a proper
congruence-preserving extension. A~positive solution to this problem 
is presented
in G. Gr\"atzer and F. Wehrung \cite{M3L}. The construction used there, the
\emph{Boolean triple construction}, as well as its generalization called the
\emph{box product construction}, see G.~Gr\"atzer and F. Wehrung 
\cite{BoxProd},
turned out to be very useful. In Section~\ref{S:BoolTr}, we discuss some of the
results that can be obtained with these constructions.

\item\textbf{Extending partial lattices to lattices}
(Sections \ref{S:Conc}, \ref{S:LiftPL}, and \ref{S:PartLat}). The
original Gr\"atzer-Schmidt solution to the characterization problem of
algebraic lattices as congruence lattices of algebras, see G. Gr\"atzer
and E.\,T. Schmidt \cite{GrSc63}, starts with a partial algebra that is further
extended to a total algebra. However, this method requires to add
infinitely many operations, thus it is,
\emph{a priori}, not suited for dealing with a class of objects of a fixed
type such as lattices. However, there are some special methods
that significantly extend the known positive results to wider classes
of distributive \jzs s.

\item\textbf{Ring-theoretical methods} (Section~\ref{S:Ring1}). A survey of the
connections between ring-theoretical problems and results and 
congruence lattice
representation problems is presented in K.\,R. Goodearl and F. 
Wehrung \cite{GoWe}.
In Section~\ref{S:Ring1}, we present a very short overview of the 
subject, as well
as a few recent results.

\item\textbf{Finitely generated varieties of lattices}
(Section~\ref{S:AlgLatt}). It is proved in M.~Plo\v{s}\v{c}ica, J. 
T\r{u}ma, and
F. Wehrung \cite{PTW} that for any nondistributive variety $\VV$ of 
lattices and any
set $X$ with at least $\aleph_2$ elements, the congruence lattice of the free
lattice $\Fg_{\VV}(X)$ satisfies many negative properties with 
respect to \CLP, see
Theorem~\ref{T:PermConURP}; in particular, its semilattice of compact 
elements is
not representable \emph{via} Schmidt's Lemma (see E.\,T. Schmidt \cite{Schm68},
Proposition~\ref{P:BasicWD}, and Theorem~\ref{T:Schmidt}) and it is 
not isomorphic
to $\Conc L$, for any sectionally complemented lattice $L$. A variety 
of lattices is
nondistributive if{f} it contains as an element either the diamond
$M_3$ or the pentagon $N_5$. As a surprising consequence, even the 
very ``simple''
finitely generated lattice varieties $\MM_3$ and $\mathbf{N}_5$ have 
complicated
congruence classes (see Definition~\ref{D:ConCl}), not completely 
understood yet.
Nevertheless, M. Plo\v{s}\v{c}ica's work is an important step in this 
direction.
\end{itemize}

\section{Uniform Refinement Properties}\label{S:URP}

The key to all known negative congruence lattice representation
results lies in considering certain infinitary axioms of the theory of
\js s (we do not need the zero in their formulation) that we call
\emph{uniform refinement properties}. Although there is no
precise definition of what a `uniform refinement property' should be in
general, the few of them that we shall review in this section 
undoubtedly offer a
very recognizable pattern.

The first idea of this pattern can be found in the case of a 
\emph{finite number}
of equations (in a given join-semilattice), of the form
  \[
  \Sigma\colon\xa_i\vee\xb_i=\mathrm{constant}
  \qquad(\mbox{for all }i\in I).
  \]
When $I=\set{i,j}$, a \emph{refinement} of $\Sigma$ can be defined as 
a collection
of four elements $\xc_{ij}^{uv}$ (for $u$, $v<2$) satisfying the equations
  \begin{equation}
  \begin{split}\label{Eq:Refab}
  \xa_i=\xc_{ij}^{00}\vee\xc_{ij}^{01}\qquad&\text{and}\qquad
  \xb_i=\xc_{ij}^{10}\vee\xc_{ij}^{11},\\
  \xa_j=\xc_{ij}^{00}\vee\xc_{ij}^{10}\qquad&\text{and}\qquad
  \xb_j=\xc_{ij}^{01}\vee\xc_{ij}^{11},
  \end{split}
  \end{equation}
see Figure~1.
\begin{figure}[t!]
  \centerline{\includegraphics{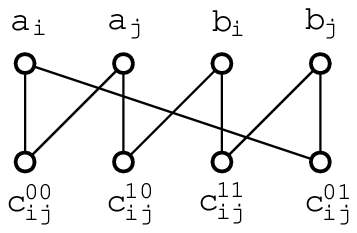}}
  \caption{Refining $\xa_i\vee\xb_i=\xa_j\vee\xb_j$.}
  \end{figure}
Observe that (\ref{Eq:Refab}) implies immediately the following:
  \begin{equation}\label{Eq:RefConseq}
  \xa_i\leq\xa_j\vee\xc_{ij}^{01}.
  \end{equation}
When $I$ is an arbitrary finite set, with powerset $\PP(I)$, one can 
extend this in
a natural way and thus define a refinement of $\Sigma$ to be a $\PP(I)$-indexed
family of elements of $S$ satisfying suitable generalizations of
\eqref{Eq:Refab}. Nevertheless, this cannot be extended immediately to the
infinite case, so that we shall focus instead on the consequence
\eqref{Eq:RefConseq} of refinement, together with
an additional ``coherence condition"
$\xc_{ik}^{01}\leq\xc_{ij}^{01}\vee\xc_{jk}^{01}$.
Thus we obtain the first uniform refinement property, see F.
Wehrung~\cite{WeURP}. This condition is a byproduct of a related 
infinitary axiom
for \emph{dimension groups} obtained in F. Wehrung \cite{We3}. 
Another (easily seen
to be equivalent) form of this axiom, denoted by $\URP_1$ in J.~T\r{u}ma and F.
Wehrung \cite{TuWe1}, is the following, we shall denote it here by $\URP$:

\begin{definition}\label{D:URP}
Let $S$ be a \js, let $\xe\in S$. We say that $S$ satisfies
\emph{$\URP$ at $\xe$}, if for all families $\famm{\xa_i}{i\in I}$ and
$\famm{\xb_i}{i\in I}$ of elements of $S$ such that
$\xa_i\vee\xb_i=\xe$, for all $i\in I$, there are elements
$\xa^*_i$, $\xb^*_i$, $\xc_{i,j}$ (for $i$, $j\in I$)
of $S$ such that the following statements hold:
\begin{enumerate}
\item $\xa^*_i\leq\xa_i$, $\xb^*_i\leq\xb_i$,
and $\xa^*_i\vee\xb^*_i=\xe$, for all $i\in I$;

\item $\xc_{i,j}\leq\xa^*_i$ and $\xc_{i,j}\leq\xb^*_j$, for all $i$,
$j\in I$;

\item $\xa^*_i\leq\xa^*_j\vee\xc_{i,j}$
and $\xb^*_j\leq\xb^*_i\vee\xc_{i,j}$, for all $i$, $j\in I$;

\item $\xc_{i,k}\leq\xc_{i,j}\vee\xc_{j,k}$, for all $i$, $j$, $k\in I$.
\end{enumerate}
We say that $S$ has $\URP$, if it has $\URP$ at all its elements.
\end{definition}

A slight weakening of $\URP$ is introduced in M. Plo\v{s}\v{c}ica,
J. T\r{u}ma, and F. Wehrung \cite{PTW}:

\begin{definition}\label{D:WURP}
Let $S$ be a \js, let $\xe\in S$. We say that $S$ satisfies
\emph{$\WURP$ at $\xe$}, if for all families $\famm{\xa_i}{i\in I}$ and
$\famm{\xb_i}{i\in I}$ of elements of $S$ such that
$\xa_i\vee\xb_i=\xe$, for all $i\in I$, there are elements
$\xc_{i,j}$ (for $i$, $j\in I$)
of $S$ such that the following statements hold:
\begin{enumerate}
\item $\xc_{i,j}\leq\xa_i,\xb_j$, for all $i$, $j\in I$;

\item $\xa_j\vee\xb_i\vee\xc_{i,j}=\xe$, for all $i$, $j\in I$;

\item $\xc_{i,k}\leq\xc_{i,j}\vee\xc_{j,k}$, for all $i$, $j$, $k\in I$.
\end{enumerate}
We say that $S$ has $\WURP$, if it has $\WURP$ at all its elements.
\end{definition}

\begin{definition}\label{D:WDhom}
Let $S$ and $T$ be \js s, let $\xe\in S$. A \jh\ $\mu\colon S\to T$
is \emph{weakly distributive}, if for all $\xa$, $\xb\in T$ such that
$\mu(\xe)=\xa\vee\xb$, there are $\xa'$, $\xb'\in S$ such that
$\mu(\xa')\leq\xa$, $\mu(\xb')\leq\xb$, and $\xa'\vee\xb'=\xe$.
\end{definition}

For further use (see Theorem~\ref{T:Schmidt}), we record here the following
definition:

\begin{definition}\label{D:DistrCon}
Let $\xa$ be a congruence of a join-semilattice $S$. We say that $\xa$ is
\begin{itemize}
\item[---] \emph{weakly distributive}, if the canonical projection 
from $S$ onto
$S/{\xa}$ is weakly distributive;

\item[---] \emph{monomial}, if every congruence class of $\xa$ has a largest
element;

\item[---] \emph{distributive}, if $\xa$ is a union of a family of weakly
distributive and monomial congruences of $S$.
\end{itemize}
For join-semilattices $S$ and $T$, a surjective \jh\ $\mu\colon S\onto T$ is
\emph{distributive}, if its kernel is a distributive congruence of $S$.
\end{definition}

The following easy result records standard facts about uniform refinement
properties and weakly distributive homomorphisms, see M. 
Plo\v{s}\v{c}ica, J. T\r{u}ma,
and F. Weh\-rung~\cite{PTW} and F. Wehrung~\cite{WeURP}:

\begin{proposition}\label{P:BasicWD}
Let $S$ and $T$ be \js s. Then the following statements hold:
\begin{enumerate}
\item If $S$ is distributive, then the set of all elements of $S$ at
which $\URP$ \pup{resp., $\WURP$} holds is closed under finite joins.

\item For any weakly distributive \jh\ $\mu\colon S\to T$ \pup{see
Definition~\textup{\ref{D:WDhom}}} and any $\xe\in S$, if $S$ has $\URP$
\pup{resp., $\WURP$} at $\xe$, then $T$ has $\URP$ \pup{resp., $\WURP$} at
$\mu(\xe)$.

\item $\URP$ implies $\WURP$.
\end{enumerate}
\end{proposition}

The following result, obtained in M. Plo\v{s}\v{c}ica, J. T\r{u}ma,
and F. Wehrung \cite{PTW} and in J. T\r{u}ma
and F. Wehrung \cite{TuWe1}, explains how uniform refinement
properties can be used to separate classes of semilattices:

\begin{theorem}\label{T:PermConURP}\hfill
\begin{enumerate}
\item For every lattice $L$ with permutable congruences, the
congruence semilattice $\Conc L$ satisfies $\URP$.

\item Let $\VV$ be a nondistributive variety of lattices, let $X$ be a
set with at least~$\aleph_2$ elements. Then $\Conc\Fg_{\VV}(X)$
does not satisfy $\WURP$.
\end{enumerate}
\end{theorem}

Although Theorem~\ref{T:PermConURP} is a difficult result, its set-theoretical
part, that explains what is so special about the cardinality $\aleph_2$, is a
very simple statement of infinite combinatorics, see C. Kuratowski 
\cite{Kura51}:

\begin{all}{The Kuratowski Free Set Theorem}
Let $n$ be a positive integer, let $X$ be a set. Then
$|X|\geq\aleph_n$ if{f} for every map $\Phi\colon[X]^n\to\fin{X}$,
there exists $U\in[X]^{n+1}$ such that 
$u\notin\Phi(U\setminus\set{u})$, for any
$u\in U$.
\end{all}

In fact, only the case $n=2$ is used.

A more complicated weakening of $\URP$ is used in J. T\r{u}ma and F.
Wehrung \cite{TuWe1}:

\begin{definition}\label{D:URPm}
Let $S$ be a \js, let $\xe\in S$. We say that $S$ satisfies
\emph{$\URP^-$ at $\xe$}, if for all families $\famm{\xa_i}{i\in I}$ and
$\famm{\xb_i}{i\in I}$ of elements of $S$ such that
$\xa_i\vee\xb_i=\xe$, for all $i\in I$, there are a subset $X$ of $I$ and
elements $\xa^*_i$, $\xb^*_i$, $\xc_{i,j}$ (for $i$, $j\in I$)
of $S$ such that the following statements hold:
\begin{enumerate}
\item $\xa^*_i\leq\xa_i$, $\xb^*_i\leq\xb_i$,
and $\xa^*_i\vee\xb^*_i=\xe$, for all $i\in I$;

\item $\xc_{i,j}\leq\xa^*_i$ and $\xc_{i,j}\leq\xb^*_j$, for all $i$,
$j\in I$;

\item $\xa^*_i\leq\xa^*_j\vee\xc_{i,j}$
and $\xb^*_j\leq\xb^*_i\vee\xc_{i,j}$, for all $i$, $j\in I$;

\item $\xc_{i,k}\leq\xc_{i,j}\vee\xc_{j,k}$, for all $i$, $j$, $k\in I$
such that the following two conditions hold:
  \begin{align*}
  i,\,k\in X&\text{ implies that }j\in X,\\
  i,\,k\notin X&\text{ implies that }j\notin X.
  \end{align*}
\end{enumerate}
We say that $S$ has $\URP^-$, if it has $\URP^-$ at all its elements.
\end{definition}

The following difficult result, see J. T\r{u}ma and F. 
Wehrung~\cite{TuWe1}, extends
Theorem~\ref{T:PermConURP}:

\begin{theorem}\label{T:PermConURPm}\hfill
\begin{enumerate}
\item For every lattice $L$ with almost permutable congruences, the
congruence semilattice $\Conc L$ satisfies $\URP^-$ at every
principal congruence of $L$.

\item Let $\VV$ be a nondistributive variety of lattices, let $X$ be a
set with at least $\aleph_2$ elements. Then the congruence
semilattice of the $\VV$-free \emph{bounded} lattice
$\Fb_{\VV}(X)$ over $X$ does
not satisfy $\URP^-$ at the largest congruence of
$\Fb_{\VV}(X)$.
\end{enumerate}
\end{theorem}

We end this section by presenting the following \emph{strengthening}
of $\URP$:

\begin{definition}\label{D:URPp}
Let $S$ be a \js, let $\xe\in S$. We say that $S$ satisfies
\emph{$\URP^+$ at $\xe$}, if for all families $\famm{\xa_i}{i\in I}$ and
$\famm{\xb_i}{i\in I}$ of elements of $S$ such that
$\xa_i\vee\xb_i=\xe$, for all $i\in I$, there are
elements $\xc_{i,j}$ (for $i$, $j\in I$)
of $S$ such that the following statements hold:
\begin{enumerate}
\item $\xc_{i,j}\leq\xa_i$ and $\xc_{i,j}\leq\xb_j$, for all $i$,
$j\in I$;

\item $\xa_i\leq\xa_j\vee\xc_{i,j}$
and $\xb_j\leq\xb_i\vee\xc_{i,j}$, for all $i$, $j\in I$;

\item $\xc_{i,k}\leq\xc_{i,j}\vee\xc_{j,k}$, for all $i$, $j$, $k\in I$.
\end{enumerate}
We say that $S$ has $\URP^+$, if it has $\URP^+$ at all its elements.
\end{definition}

In view of Corollary~\ref{C:Informal}, it follows that all known recent
representation theorems yield semilattices with $\URP^+$:

\begin{proposition}\label{P:RelCplURP+}
The congruence semilattice $\Conc L$ satisfies $\URP^+$, for any
relatively complemented lattice $L$.
\end{proposition}

\begin{proof}
Let $\xe\in\Conc L$, let $\famm{\xa_i}{i\in I}$ and
$\famm{\xb_i}{i\in I}$ be families of elements of $\Conc L$ such that
$\xa_i\vee\xb_i=\xe$, for all $i\in I$, we shall find elements
$\xc_{i,j}$ (for all $i$, $j\in I$) of $\Conc L$ that satisfy the required
inequalities. Since $L$ is relatively complemented, every compact congruence
of $L$ is principal, thus there are elements $u\leq v$ in $L$ such that
$\xe=\Theta(u,v)$.

Furthermore, it follows from Proposition~3.2 in F. 
Wehrung~\cite{WeURP} that $L$
is ``congruence splitting'', thus, for all $i\in I$, since 
$\xa_i\vee\xb_i=\xe$,
there are $x_i$, $y_i\in[u,v]$ such that $x_i\vee y_i=v$,
$\Theta(u,x_i)\subseteq\xa_i$, and $\Theta(u,y_i)\subseteq\xb_i$.

Since $\xa_i$ is a compact congruence of $L$, it is principal, thus
we can write $\xa_i=\Theta(u_i,v_i)$,
for some $u_i\leq v_i$ in $L$, so 
$\Theta(u_i,v_i)\subseteq\Theta(u,v)$. Thus there
exists a subdivision of the interval $[u_i,v_i]$ whose subintervals all
weakly project into $[u,v]$ (see Theorem~III.1.2 in G. Gr\"atzer \cite{GLT2}).
Hence, since $L$ is relatively complemented, any of these intervals 
is projective
to a subinterval of $[u,v]$ (see Exercise~III.1.3 in G. Gr\"atzer 
\cite{GLT2}), thus
(again since $L$ is relatively complemented) to an interval of the 
form $[u,w]$,
where $u\leq w\leq v$. Denoting by $s_i$ the join of all the $w$-s thus
obtained; we get that $u\leq s_i\leq v$, while
$\xa_i=\Theta(u,s_i)$. Similarly, we can get $t_i\in[u,v]$ such that
$\xb_i=\Theta(u,t_i)$. Define $a_i=x_i\vee s_i$ and $b_i=y_i\vee 
t_i$. The relevant
properties of $a_i$ and $b_i$ are the following:
  \begin{equation}\label{Eq:Relaibi}
  a_i,b_i\in[u,v];\qquad a_i\vee b_i=v;
  \qquad\xa_i=\Theta(u,a_i);\qquad\xb_i=\Theta(u,b_i).
  \end{equation}
We define compact congruences of $L$ by
  \[
  \xc'_{i,j}=\Theta^+(a_i,a_j)\qquad\text{and}\qquad\xc''_{i,j}=\Theta^ 
+(b_j,b_i),
  \qquad\text{for all }i,\,j\in I,
  \]
where we define $\Theta^+(x,y)=\Theta(x\wedge y,x)$, see 
Section~\ref{S:Conc}. For
all $i$, $j$, $k\in I$, it is not hard to verify the following inequalities:
  \begin{align*}
  \xc'_{i,j}&\subseteq\xa_i,\,\xb_j;&\xa_i&\subseteq\xa_j\vee\xc'_{i,j};&
  \xc'_{i,k}&\subseteq\xc'_{i,j}\vee\xc'_{j,k};\\
  \xc''_{i,j}&\subseteq\xa_i,\,\xb_j;&\xb_j&\subseteq\xb_i\vee\xc''_{i,j};&
  \xc''_{i,k}&\subseteq\xc''_{i,j}\vee\xc''_{j,k}.
  \end{align*}
For example, for any congruence $\xe$ of $L$, if $\xb_j\subseteq\xe$, that is,
$b_j\equiv_{\xe}u$, then, by the first two equations of \eqref{Eq:Relaibi},
$a_j\equiv_{\xe}v$, whence $a_i\leq_{\xe}a_j$, that is,
$\xc'_{i,j}\subseteq\xb_j$.
Therefore, by putting $\xc_{i,j}=\xc'_{i,j}\vee\xc''_{i,j}$, for all 
$i$, $j\in I$,
we obtain the following inequalities
  \[
  \xc_{i,j}\subseteq\xa_i,\xb_j;\qquad\xa_i\subseteq\xa_j\vee\xc_{i,j};\qquad
  \xb_j\subseteq\xb_i\vee\xc_{i,j};\qquad\xc_{i,k}\subseteq\xc_{i,j}\vee\xc_{j,k},
  \]
which concludes the proof.
\end{proof}

We also refer the reader to Problem~\ref{Pb:URPsc}.

We can also prove the following result (compare it with
Theorem~\ref{T:New0Dim}(i)):

\begin{proposition}\label{P:WOlim}
Any $\vee$-direct limit over a \emph{totally ordered set} of
distributive \emph{lattices} satisfies $\URP^+$.
\end{proposition}

\begin{proof}
Let $S$ be a $\vee$-direct limit over a totally ordered set
$\theta$ of distributive lattices, say,
$S=\varinjlim_{\alpha\in\theta}S_\alpha$ of distributive lattices
$S_\alpha$, with transition maps
$f_{\alpha,\beta}\colon S_\alpha\to S_\beta$, for
$\alpha\leq\beta$ in $\theta$ and limiting maps
$f_\alpha\colon S_\alpha\to S$, for $\alpha\in\theta$ (so the
$f_{\alpha,\beta}$-s and the $f_\alpha$-s are \jh s). We prove that
$S$ satisfies $\URP^+$.
By extracting from $\theta$ a cofinal well-ordered chain, we may
assume, without loss of generality, that $\theta$ is an
\emph{ordinal}. Let $\xe\in S$, let $\famm{\xa_i}{i\in I}$ and
$\famm{\xb_i}{i\in I}$ be families of elements of $S$ such that
$\xa_i\vee\xb_i=\xe$, for all $i\in I$; we shall find elements
$\xc_{i,j}$ (for $i$, $j\in I$) of $S$ that satisfy the required
inequalities.

Without loss of generality, $\xe$ belongs to the range of $f_0$, so
$\xe=f_0(\xe^0)$, for some $\xe^0\in S_0$. Put
$\xe^\alpha=f_{0,\alpha}(\xe^0)$, for all $\alpha<\theta$; observe that
$f_\alpha(\xe^\alpha)=\xe$.

For all $i\in I$, there are $\alpha<\theta$ and $\xu$,
$\xv\in S_\alpha$ such that $\xa_i=f_\alpha(\xu)$,
$\xb_i=f_\alpha(\xv)$, and $\xu\vee\xv=\xe^\alpha$. Denote by
$\mu(i)$ the least such $\alpha$, and let
$\seq{\xa_i^{\mu(i)},\xb_i^{\mu(i)}}$ be a corresponding choice for
$\seq{\xu,\xv}$. Put $\xa_i^\alpha=f_{\mu(i),\alpha}(\xa_i^{\mu(i)})$
and $\xb_i^\alpha=f_{\mu(i),\alpha}(\xb_i^{\mu(i)})$, for all
$\alpha<\theta$ with $\alpha\geq\mu(i)$. Further, define
$\nu(i,j)=\max\set{\mu(i),\mu(j)}$, for all $i$, $j\in I$.

We denote by $\wedge_\alpha$ the meet operation in $S_\alpha$, for all
$\alpha<\theta$, and for all $i$, $j\in I$, we define
  \[
  \xc_{i,j}=f_{\nu(i,j)}
  \bigl(\xa_i^{\nu(i,j)}\wedge_{\nu(i,j)}\xb_j^{\nu(i,j)}\bigr).
  \]
To conclude the proof, if suffices to establish that the elements
$\xc_{i,j}$ thus defined satisfy the required inequalities.

Let $i$, $j$, $k\in I$, set $\alpha=\nu(i,j)$, $\beta=\nu(j,k)$, and
$\gamma=\nu(i,k)$. We observe that the following inequality holds:
  \begin{equation}\label{Eq:triang}
  \gamma\leq\max\set{\alpha,\beta}.
  \end{equation}
We first observe that $\xc_{i,j}\leq f_\alpha(\xa_i^\alpha)=\xa_i$,
and, similarly, $\xc_{i,j}\leq\xb_j$.

We further compute:
  \begin{align*}
  \xa_i^\alpha&=\xa_i^\alpha\wedge_\alpha\xe^\alpha\\
  &=\xa_i^\alpha\wedge_\alpha(\xa_j^\alpha\vee\xb_j^\alpha)\\
  &=(\xa_i^\alpha\wedge_\alpha\xa_j^\alpha)\vee
  (\xa_i^\alpha\wedge_\alpha\xb_j^\alpha)&&
  (\text{by the distributivity of }S_\alpha)\\
  &\leq\xa_j^\alpha\vee(\xa_i^\alpha\wedge_\alpha\xb_j^\alpha),
  \end{align*}
whence, by applying $f_\alpha$, we obtain that
$\xa_i\leq\xa_j\vee\xc_{i,j}$. The proof of
$\xb_j\leq\xb_i\vee\xc_{i,j}$ is similar.

Finally, we verify the inequality
  \begin{equation}\label{Eq:Trcij}
  \xc_{i,k}\leq\xc_{i,j}\vee\xc_{j,k}.
  \end{equation}
We separate cases.
\smallskip

\noindent\textbf{Case 1.} $\alpha\leq\beta$.
It follows from \eqref{Eq:triang} that $\gamma\leq\beta$ as well. We
establish further inequalities. We begin with the following:
  \begin{equation}\label{Eq:Trcij11}
  f_{\gamma,\beta}(\xa_i^\gamma\wedge_\gamma\xb_k^\gamma)\leq
  f_{\alpha,\beta}(\xa_i^\alpha\wedge_\alpha\xb_j^\alpha)
  \vee\xa_j^\beta.
  \end{equation}
Indeed,
$f_{\gamma,\beta}(\xa_i^\gamma\wedge_\gamma\xb_k^\gamma)\leq
f_{\gamma,\beta}(\xa_i^\gamma)=\xa_i^\beta=
f_{\alpha,\beta}(\xa_i^\alpha)$ and
$\xa_j^\beta=f_{\alpha,\beta}(\xa_j^\alpha)$, thus, in order to prove
\eqref{Eq:Trcij11}, it suffices to verify that
$\xa_i^\alpha\leq
(\xa_i^\alpha\wedge_\alpha\xb_j^\alpha)\vee\xa_j^\alpha$, which holds
by the distributivity of $S_\alpha$ since
$\xa_j^\alpha\vee\xb_j^\alpha=\xe^\alpha$.

Next, we prove the following inequality:
  \begin{equation}\label{Eq:Trcij12}
  f_{\gamma,\beta}(\xa_i^\gamma\wedge_\gamma\xb_k^\gamma)\leq
  f_{\alpha,\beta}(\xa_i^\alpha\wedge_\alpha\xb_j^\alpha)
  \vee\xb_k^\beta.
  \end{equation}
Indeed,
$f_{\gamma,\beta}(\xa_i^\gamma\wedge_\gamma\xb_k^\gamma)\leq
f_{\gamma,\beta}(\xb_k^\gamma)=\xb_k^\beta\leq
f_{\alpha,\beta}(\xa_i^\alpha\wedge_\alpha\xb_j^\alpha)\vee\xb_k^\beta$.
Therefore, by using the distributivity of $S_\beta$ and the inequalities
\eqref{Eq:Trcij11} and \eqref{Eq:Trcij12}, we obtain the following 
inequalities:
  \begin{align*}
  \xc_{i,j}\vee\xc_{j,k}&=
  f_\alpha\bigl(\xa_i^\alpha\wedge_\alpha\xb_j^\alpha\bigr)\vee
  f_\beta\bigl(\xa_j^\beta\wedge_\beta\xb_k^\beta\bigr)\\
  &=f_\beta\Bigl(
  f_{\alpha,\beta}\bigl(\xa_i^\alpha\wedge_\alpha\xb_j^\alpha\bigr)
  \vee
  \bigl(\xa_j^\beta\wedge_\beta\xb_k^\beta\bigr)\Bigr)\\
  &=f_\beta\Bigl(\bigl(
  f_{\alpha,\beta}\bigl(\xa_i^\alpha\wedge_\alpha\xb_j^\alpha\bigr)
  \vee\xa_j^\beta\bigr)\wedge_\beta
  \bigl(f_{\alpha,\beta}\bigl(\xa_i^\alpha\wedge_\alpha\xb_j^\alpha\bigr)
  \vee\xb_k^\beta\bigr)\Bigr)\\
  &\geq f_\beta\bigl(f_{\gamma,\beta}\bigl(
  \xa_i^\gamma\wedge_\gamma\xb_k^\gamma\bigr)\bigr)\\
  &=f_\gamma\bigl(\xa_i^\gamma\wedge_\gamma\xb_k^\gamma\bigr)\\
  &=\xc_{i,k},
  \end{align*}
thus obtaining \eqref{Eq:Trcij}.

\smallskip
\noindent\textbf{Case 2.} $\beta\leq\alpha$.
It follows from \eqref{Eq:triang} that $\gamma\leq\alpha$ as well.
In a fashion similar to Case 1, one can prove the following inequalities
  \begin{align*}
  f_{\gamma,\alpha}\bigl(\xa_i^\gamma\wedge_\gamma\xb_k^\gamma\bigr)
  &\leq\xb_j^\alpha\vee
  f_{\beta,\alpha}\bigl(\xa_j^\beta\wedge_\beta\xb_k^\beta\bigr),\\
  f_{\gamma,\alpha}\bigl(\xa_i^\gamma\wedge_\gamma\xb_k^\gamma\bigr)
  &\leq\xa_i^\alpha\vee
  f_{\beta,\alpha}\bigl(\xa_j^\beta\wedge_\beta\xb_k^\beta\bigr),
  \end{align*}
thus, as before, obtaining \eqref{Eq:Trcij}.
\end{proof}

We summarize in Table~1 many known results and questions about uniform
refinement properties, sometimes anticipating some
subsequent sections of the present paper. We use the following abbreviations:
\begin{itemize}
\item distr. image of gBs = distributive image of a
generalized Boolean semilattice;

\item $L$ p.c. = $L$ with permutable congruences;

\item $L$ s.c. = $L$ sectionally complemented;

\item $L$ a.p.c. = $L$ with almost permutable congruences;

\item $L$ r.c. = $L$ relatively complemented;

\item $F_{\mathrm{b}}=\Fb_{\VV}(\omega_2)$, where $\VV$ is a nondistributive
variety of lattices;

\item $F=\Fg_{\VV}(\omega_2)$, where $\VV$ is a nondistributive
variety of lattices.
\end{itemize}

Also, the entry of the table marked by $^{(*)}$ means that $\URP^-$ holds in
$\Conc L$ at principal congruences of $L$. Finally, we recall that 
$\URP^+$ implies
$\URP$, which implies both $\WURP$ and $\URP^-$.
  \begin{table}[hbt]
  \begin{tabular}{|l|c|c|c|c|c|c|c|}
  \cline{2-8}
  \multicolumn{1}{l|}{} & $\Conc F_{\mathrm{b}}$ & $\Conc F$ &
  $\vcenter{\hbox{\strut distr. image}\hbox{of gBs\strut}}$ &
  $\vcenter{\hbox{\strut$\Conc L$,}\hbox{$L$ p.c.\strut}}$ &
  $\vcenter{\hbox{\strut$\Conc L$,}\hbox{$L$ s.c.\strut}}$ &
  $\vcenter{\hbox{\strut$\Conc L$,}\hbox{$L$ a.p.c.\strut}}$ &
  $\vcenter{\hbox{\strut$\Conc L$,}\hbox{$L$ r.c.\strut}}$\tvi\\
  \hline
  $\URP^+$\tvi & No & No & ? & ? & ? & ? & Yes\\
  \hline
  $\URP$\tvi & No & No & Yes & Yes & Yes & ? & Yes\\
  \hline
  $\WURP$\tvi & No & No & Yes & Yes & Yes & ? & Yes\\
  \hline
  $\URP^-$\tvi & No & No & Yes & Yes & Yes & Yes$^{(*)}$ & Yes\\
  \hline
  \end{tabular}
  \caption{\tvi Uniform refinement properties and congruence semilattices.}
  \end{table}

\section{The $\bt{L}$ construction, tensor product, and box product}
\label{S:BoolTr}

For a lattice $L$, we define $\bt{L}$ as the set of all triples
$\seq{x,y,z}\in L^3$ that are \emph{Boolean}, that is, the following equalities
hold:
  \begin{align*}
    x  &= (x \vee y) \wedge  (x \vee z),\\
    y  &= (y \vee x) \wedge  (y \vee z),\\
    z  &= (z \vee x) \wedge  (z \vee y).
  \end{align*}
The set $\bt{L}$ is endowed with the restriction of the componentwise ordering
on~$L^3$. It can be shown that $\bt{L}$ is a closure system in $L^3$, 
thus it is a
lattice. Some of the relevant information about this construction is 
summarized in
the following result, see G. Gr\"atzer and F. Wehrung \cite{M3L}.

\begin{proposition}\label{P:basic}
For any lattice $L$, the following statements hold:
\begin{enumerate}
\item For a congruence $\xa$ of $L$, let $\xa^3$ denote the
congruence of $L^3$ defined as $\xa$ componentwise. Let $\bt{\xa}$ be the
restriction of $\xa^3$ to $\bt{L}$. Then $\bt{\xa}$ is a congruence of
$\bt{L}$, and every congruence of $\bt{L}$ is of the form $\bt{\xa}$, for a
unique congruence $\xa$ of $L$.

\item The diagonal map $x\mapsto\seq{x,x,x}$ is a 
congruence-preserving embedding
from~$L$ into $\bt{L}$. If, in addition, $L$ has a zero \pup{resp., a 
unit}, then
the map $x\mapsto\seq{x,0,0}$ \pup{resp., $x\mapsto\seq{1,x,x}$} is a
congruence-preserving embedding from $L$ into $\bt{L}$ whose range is an ideal
\pup{resp., a dual ideal} of~$L$.
\end{enumerate}
\end{proposition}

This solves the question raised above, namely, if $L$ is a nontrivial lattice,
then the diagonal map from $L$ into $\bt{L}$ defines a proper
congruence-preserving extension of $L$.

We say that a lattice $L$ is \emph{regular}, if any two congruences of $L$ that
share a congruence class are equal. By iterating the $\bt{L}$ construction, its
refinement $\bt{L,a}$ (the latter is a convex sublattice of $\bt{L}$), and the
gluing construction, G.~Gr\"atzer and E.\,T. Schmidt prove in~\cite{RegLat} the
following theorem:

\begin{theorem}\label{T:general}
Every lattice $L$ has a congruence-preserving embedding into a regular lattice
$\tilde{L}$. If $L$ has a zero, then one can suppose that $\tilde{L}$ 
has a zero
and that $0_L=0_{\tilde{L}}$.
\end{theorem}

It is also observed in the same paper that every compact congruence 
of a regular
lattice is principal. Hence, \emph{if \CLP\ can be solved positively, 
then it can
be solved with lattices in which every compact congruence is principal}.

The $\bt{L}$ construction has a far reaching generalization, the 
\emph{box product}
of lattices. The box product improves the classical \emph{tensor product} of
join-semilattices with zero, see, for example, G. Gr\"atzer, H. 
Lakser, and R.\,W.
Quackenbush \cite{GLQ81} and G. Gr\"atzer and F. Wehrung 
\cite{TensRev}. For \jzs s
$A$ and $B$, the tensor product $A\otimes B$ is defined in a fashion 
similar to the
tensor product of vector spaces in linear algebra, in particular, it is also a
\jzs. However, even in case both $A$ and $B$ are lattices, $A\otimes B$ is not
necessarily a lattice, see G.~Gr\"atzer and F. Wehrung \cite{Transf,M3Dnmod}.

For a lattice $L$, we put $\bot_L=\set{0}$, if $L$ has a zero (least 
element) $0$,
and $\bot_L=\es$, otherwise. For lattices $A$ and $B$ and
$\seq{a,b}\in A\times B$, we define
  \begin{align*}
  \bot_{A,B}&=(A\times\bot_B)\cup(\bot_A\times B),\\
  a\ltp b&=\bot_{A,B}\cup\setm{\seq{x,y}\in A\times B}{x\leq a\text{ 
and }y\leq b},\\
  a\bp b&=\setm{\seq{x,y}\in A\times B}{x\leq a\text{ or }y\leq b}.
  \end{align*}
We denote by $A\bp B$ the \emph{box product} of $A$ and $B$;
the elements of $A\bp B$ are subsets of $A\times B$ that can be represented as
finite intersections of sets of the form $a\bp b$, with $a \in A$ and 
$b \in B$.
Unlike the tensor product $A\otimes B$, it is always a lattice, see 
G. Gr\"atzer and
F. Wehrung \cite{BoxProd}. An element of $A\bp B$ is \emph{confined}, if it is
contained in some element of the form $a\ltp b$, for $\seq{a,b}\in 
A\times B$. The
ideal $A\ltp B$ of all confined elements of $A\bp B$ is nonempty 
if{f} either $A$ or
$B$ is bounded, or both $A$ and $B$ have a zero, or both $A$ and $B$ 
have a unit;
we call it the \emph{lattice tensor product} of $A$ and~$B$, and then the
\emph{Isomorphism Theorem} holds, which implies the following formula:
  \begin{equation}\label{Eq:IsFla}
  (\Conc A)\otimes(\Conc B)\cong\Conc(A\ltp B).
  \end{equation}
In particular, for a lattice $L$, the lattice tensor product $M_3\ltp L$ is
isomorphic to the lattice $\bt{L}$ introduced at the beginning of
Section~\ref{S:BoolTr}. Such an isomorphism is called a 
\emph{coordinatization},
and sometimes provides a more convenient way to compute in lattice tensor
products. Arbitrary lattice tensor products of bounded lattices are 
coordinatized
in G. Gr\"atzer and M. Greenberg \cite{Coord1,Coord2,Coord3,Coord4}.

To prove the Isomorphism Theorem, one first needs to verify
a more general result that extends the formula
\eqref{Eq:IsFla} to so-called \emph{capped sub-tensor products} of 
$A$ and $B$, this
is the main result of G. Gr\"atzer and F. Wehrung \cite{TensRev}. 
Then, one needs to
verify that if $A$ and $B$ are lattices with zero, then $A\ltp B$ is a capped
sub-tensor product of $A$ and $B$, see G. Gr\"atzer and F. Wehrung 
\cite{BoxProd}.
The isomorphism of \eqref{Eq:IsFla} carries 
$\Theta_A(a,a')\otimes\Theta_B(b,b')$ to
$\Theta_{A\ltp B}((a\ltp b')\vee(a'\ltp b),a'\ltp b')$, for all 
$a\leq a'$ in $A$ and
$b\leq b'$ in $B$.

This has an application to the following problem.
We say that a lattice $L$ is an \emph{automorphism-preserving extension} of
a sublattice $K$, if every automorphism of $K$ extends to a unique 
automorphism of
$L$ and $K$ is closed under all automorphisms of $L$.
By iterating the box product construction together with gluing, 
G.~Gr\"atzer and
F. Wehrung solve in \cite{StrIndep} a problem already proposed in the 
first edition
of the monograph G. Gr\"atzer~\cite{GLT2}, by proving the following:

\begin{theorem}[The Strong Independence Theorem for arbitrary
lattices]\label{T:StrIndep}
For every nontrivial lattice $L_{\mathrm{C}}$ and every
lattice $L_{\mathrm{A}}$, there exists a lattice $L$ that is both a
congruence-preserving extension of
$L_{\mathrm{C}}$ and an automorphism-preserving extension of
$L_{\mathrm{A}}$. Furthermore, if both $L_{\mathrm{C}}$ and 
$L_{\mathrm{A}}$ have
a zero, then $L$ can be taken a zero-preserving extension of both 
$L_{\mathrm{C}}$
and $L_{\mathrm{A}}$.
\end{theorem}

Because of the well-known result of G. Birkhoff that states that every group
appears as the automorphism group of some lattice, it follows that for every
nontrivial lattice $K$ and every group $G$, there exists a 
congruence-preserving
extension $L$ of~$K$ such that $\Aut L\cong G$. Observe that $\Con 
K\cong\Con L$.

The essential difficulty of the proof of Theorem~\ref{T:StrIndep} lies in the
construction, for a given lattice $L_{\mathrm{C}}$, of a \emph{rigid},
congruence-preserving extension $\ol{L}_{\mathrm{C}}$ of $L_{\mathrm{C}}$. This
construction is performed in several steps. For a bounded lattice 
$L$, we denote by
$M_3\lfloor L\rfloor$ the set of all triples $\seq{x,y,z}$ of 
$\bt{L}$ such that
either $x=0$ or $x=1$, partially ordered under inclusion. Then
$M_3\lfloor L\rfloor$ is a lattice, and its congruence lattice is
isomorphic to the lattice of all congruences of $L$ that are either
coarse or for which the congruence class of zero is zero.

Moreover, by using the lattice tensor product, one associates, with every
lattice~$L_{\mathrm{C}}$, a ``large enough'' simple, bounded lattice
$S$ (whose cardinality may be larger than that of $L_{\mathrm{C}}$). Put
$T=M_3\lfloor S\rfloor$. For a principal dual ideal $J$ of
$L_{\mathrm{C}}$, we glue~$L_{\mathrm{C}}$, with the dual ideal
$J$, with $V=T\ltp J$, with the ideal $p\ltp J$, where $p$ denotes the unique
atom of $T$. The result of this construction is a 
congruence-preserving extension of
$L_{\mathrm{C}}$. By iterating this construction transfinitely many 
times, we obtain
a rigid, congruence preserving-extension $\ol{L}_{\mathrm{C}}$ of 
$L_{\mathrm{C}}$.
Observe that the cardinality of $\ol{L}_{\mathrm{C}}$ may be larger than the
cardinality of $L_{\mathrm{C}}$. Furthermore, the extension 
$\ol{L}_{\mathrm{C}}$
thus constructed has a strong indecomposability property called 
\emph{steepness}.

By using much more elementary techniques introduced earlier in G. Gr\"atzer and
E.\,T. Schmidt \cite{GS95b}, for every lattice $L_{\mathrm{A}}$, one 
can construct a
simple, automorphism-preserving extension $\ol{L}_{\mathrm{A}}$ of
$L_{\mathrm{A}}$. Now, if $L_{\mathrm{A}}$ and $L_{\mathrm{C}}$ are given, the
extensions $\ol{L}_{\mathrm{A}}$ and $\ol{L}_{\mathrm{C}}$ are 
constructed, then we
put $L=\ol{L}_{\mathrm{A}}\ltp\ol{L}_{\mathrm{C}}$. Since 
$\ol{L}_{\mathrm{A}}$ is
simple, $L$ is a congruence-preserving extension of 
$\ol{L}_{\mathrm{C}}$, thus of
$L_{\mathrm{C}}$. Furthermore, every automorphism of $L_{\mathrm{A}}$ 
induces an
automorphism of $\ol{L}_{\mathrm{A}}$, thus an automorphism of $L$. 
By using the
steepness of $\ol{L}_{\mathrm{C}}$, one can prove, and this is the 
hardest part of
the proof, that there are no other automorphisms of $L$. Therefore, $L$ is an
automorphism-preserving extension of $\ol{L}_{\mathrm{A}}$.

If $X$ and $Y$ are subsets of a lattice $L$, we say that a map
$\varphi\colon X\to Y$ is \emph{algebraic}, if there exists a lattice
polynomial $\mathbf{p}$, with one variable and with parameters 
from~$L$, such that
$\varphi(x)=\mathbf{p}(x)$, for all $x\in X$. In G. Gr\"atzer and E.\,T.
Schmidt \cite{NewAppr1}, the following result is established:

\begin{theorem}\label{T:IntEquiv}
Let $K$ be a bounded lattice, let $[a,b]$ and $[c,d]$ be intervals of $K$, and
let $\varphi\colon[a,b]\to[c,d]$ be an isomorphism between these two intervals.
Then $K$ has a $\seq{\vee,\wedge,\varphi,\varphi^{-1}}$-congruence-preserving
extension into a bounded lattice $L$ such that both $\varphi$ and 
$\varphi^{-1}$
are algebraic in $L$, and $K$ is a convex sublattice of $L$.
In particular, the congruence lattice of the partial algebra
$\seq{K,\vee,\wedge,\varphi,\varphi^{-1}}$ is isomorphic to the 
congruence lattice
of the bounded lattice $\seq{L,\vee,\wedge}$.
\end{theorem}

The construction of Theorem~\ref{T:IntEquiv} uses refinements of the $\bt{L}$
construction together with the box product construction and gluing. 
Furthermore, it
does not require transfinite induction, in particular, it preserves finiteness.
Starting with a relatively complemented lattice $K$, the extension 
$L$ constructed
by Theorem~\ref{T:IntEquiv} is not relatively complemented. Compare
with Corollary~\ref{C:Informal}.

This result is extended to a family of isomorphisms in G. Gr\"atzer and E.\,T.
Schmidt \cite{NewAppr1}, and to a family of surjective homomorphisms between
intervals without requiring the inverses (that is, $\varphi^{-1}$) in the
extended language, see G. Gr\"atzer, M. Greenberg, and E.\,T. Schmidt
\cite{NewAppr2}. The latter construction involves the consideration 
of the lattice
tensor product of the original lattice not with $M_3$, but with $N_6$, the six
element sectionally complemented lattice obtained by replacing one of the lower
prime intervals of the square by a square. The manipulation of the 
elements of the
box product is made more convenient by the coordinatization of
lattice tensor products studied in G. Gr\"atzer and M. Greenberg
\cite{Coord1,Coord2,Coord3,Coord4}.

However, these methods alone are not sufficient to solve \CLP, because of the
following observation. They extend a lattice $K$ to a lattice $L$ 
whose congruence
lattice is isomorphic to the lattice $\Con^{\Phi}K$ of all lattice 
congruences of
$K$ having the substitution property with respect to all the 
operations of $\Phi$,
where $\Phi$ is a set of partial unary functions on $K$. In 
particular, $\Con L$ is
isomorphic to an \emph{algebraic} subset of $\Con K$ (a subset $X$ of 
a lattice $A$
is \emph{algebraic}, if it is closed under arbitrary meets and 
nonempty directed
joins of $A$), hence $\Conc L$ is the image of $\Conc K$ under a 
weakly distributive
\jzh\ (namely, the one that with a compact congruence $\xa$ associates the
$\Phi$-congruence generated by~$\xa$), see Section~\ref{S:URP}. In 
particular, if
$K$ is already obtained from an already known representation theorem, 
then $\Conc K$
satisfies the axiom $\URP$ (see Corollary~\ref{C:Informal} and the 
comments that
follow it), thus, by Proposition~\ref{P:BasicWD}(ii), $\Conc L$ also
satisfies~$\URP$.

Still this does not rule out the following possible approach of \CLP, hinted at
in Problem~1 in G. Gr\"atzer and E.\,T. Schmidt \cite{NewAppr1}. If 
one could prove
that every algebraic distributive lattice is isomorphic to some lattice of the
form $\Con^{\Phi}K$, where $K$ is a lattice with zero and $\Phi$ is a 
set of partial
surjective homomorphisms between intervals of $K$ satisfying certain simple
conditions, then \CLP\ would be solved positively. Observe that this 
would imply
that every algebraic distributive lattice $D$ is isomorphic to an 
algebraic subset
of some algebraic lattice of the form $\Con K$; hence, if the 
semilattice of compact
elements of $D$ does not satisfy $\URP$, then neither does the semilattice
$\Conc K$. Thus a natural guess would be to start with~$K$ being, say, a free
lattice.

The Isomorphism Theorem (see \eqref{Eq:IsFla}) has another interesting
consequence, see G. Gr\"atzer and F. Wehrung \cite{BoxProd}. We
say that a \jzs\ $S$ is \emph{$\seq{0}$-rep\-re\-sent\-able} (resp.,
\emph{$\seq{0,1}$-representable}), if there exists a lattice $L$ with 
zero (resp.,
a bounded lattice $L$) such that $S\cong\Conc L$. The problem whether any
representable semilattice $S$ (i.e., a semilattice $S$ for which there exists a
lattice $L$ such that $S\cong\Conc L$) is $\seq{0}$-representable is open.

\begin{theorem}\label{T:TensProd}
Let $S$ and $T$ be \jzs s. Then the following statements hold:
\begin{enumerate}
\item If both $S$ and $T$ are $\seq{0}$-representable, then $S\otimes T$ is
$\seq{0}$-representable.

\item If both $S$ and $T$ are $\seq{0,1}$-representable, then $S\otimes T$ is
$\seq{0,1}$-representable.

\item If $S$ is representable and $T$ is $\seq{0,1}$-representable, then
$S\otimes T$ is representable.
\end{enumerate}
\end{theorem}

This result can be easily extended to \emph{iterated tensor products} of
\jzus s. For \jzus s $S$ and $T$, the rule $x\mapsto x\otimes 1_T$ defines a
\jzu-embedding from $S$ into $S\otimes T$.
For a family $\famm{S_i}{i\in I}$ of \jzus s and finite subsets 
$I_0\subseteq I_1$
of $I$, one defines similarly a \jzu-embedding from $\bigotimes_{i\in I_0}S_i$
into $\bigotimes_{i\in I_1}S_i$. These maps obviously form a direct system of
\jzus s and \jzu-embeddings; let $\bigotimes_{i\in I}S_i$ denote its 
direct limit,
the \emph{iterated tensor product} of the $S_i$-s.
Suppose now that $S_i=\Conc L_i$, for some bounded lattice $L_i$, for 
all $i\in I$.
By arguing as for semilattices except that $\otimes$ is replaced by $\ltp$, we
obtain a $\seq{0,1}$-lattice embedding from $\ltp_{i\in I_0}L_i$ into
$\ltp_{i\in I_1}L_i$. These maps also form a direct system; denote its direct
limit by $\ltp_{i\in I}L_i$. Since the $\Conc$ functor preserves direct limits
(see Proposition~\ref{P:PresLim}), we obtain the formula
  \[
  \Conc\bigl(\ltp_{i\in I}L_i\bigr)\cong\bigotimes_{i\in I}(\Conc L_i).
  \]
This yields the following result:

\begin{theorem}\label{T:ItTensProd}
Any iterated tensor product of $\seq{0,1}$-representable \jzs s is
$\seq{0,1}$-representable.
\end{theorem}

This result is similar to K.\,R. Goodearl and D.\,E. Handelman
\cite[Theorem~3.5]{GoHa}, where it is proved that for any family
$\famm{G_i}{i\in I}$ of dimension groups with order-unit, if every $G_i$ is
isomorphic to the $K_0$ of some locally matricial algebra (over a given field),
then so is the iterated tensor product $\bigotimes_{i\in I}G_i$.

\section{The functor $\Conc$ on partial lattices}
\label{S:Conc}

Throughout the paper we shall make use of the following categories:

\begin{itemize}
\item $\CL$, the category of all lattices and lattice homomorphisms;

\item $\CS$, the category of all \jzs s and \jzh s;

\item $\CSd$, the full subcategory of $\CS$ whose objects are the
\emph{distributive} \jzs s;

\item $\CSfd$, the full subcategory of $\CS$ whose objects are the
\emph{finite distributive} \jzs s;

\item $\CSfb$, the full subcategory of $\CS$ whose objects are the
\emph{finite Boolean} \jzs s.
\end{itemize}

The correspondence that with every lattice $L$ associates its
congruence semilattice $\Conc L$ can be extended to a \emph{functor}
from $\CL$ to $\CS$. For a lattice
homomorphism $f\colon K\to L$, let $\Conc f$ be the map from $\Conc K$ to
$\Conc L$ that with every compact congruence $\xa$ associates the
congruence of $L$ generated by all pairs $\seq{f(x),f(y)}$, where
$\seq{x,y}\in\xa$.

This can be easily extended to \emph{partial lattices}.
The precise concepts are summarized in the following two definitions (see
F. Wehrung \cite{Wehr1con}). We observe that the definition of a 
partial lattice
that we use here is very closely related to the one used in R. Freese, J.
Je\v{z}ek, and J.\,B. Nation \cite{FJN} but not to the one in G. Gr\"atzer
\cite{GLT2}.

\begin{definition}\label{D:PartLatt}\hfill
\begin{enumerate}
\item A \emph{partial prelattice} is a structure
$\seq{P,\leq,\bigvee,\bigwedge}$, where $P$ is a nonempty set, $\leq$
is a quasi-ordering on~$P$, and $\bigvee$, $\bigwedge$ are partial
functions from the set $\fine P$ of all nonempty finite subsets of $P$
to $P$ satisfying the following properties:
\begin{enumerate}
\item $a=\bigvee X$ implies that $a=\sup X$,
for all $a\in P$ and all $X\in\fine P$.

\item $a=\bigwedge X$ implies that $a=\inf X$,
for all $a\in P$ and all $X\in\fine P$.
\end{enumerate}

(By $a=\sup X$, we mean that an element $b$ of $P$ is an upper bound 
of $X$ if{f}
$a\leq b$. The statement $a=\inf X$ is defined dually.)

\item $P$ is a \emph{partial lattice}, if $\leq$ is antisymmetric.

\item A \emph{congruence} of~$P$ is a quasi-ordering $\preceq$ of~$P$
containing $\leq$ such that $\seq{P,\preceq,\bigvee,\bigwedge}$ is a
partial prelattice.

\end{enumerate}
\end{definition}

Lattices are naturally identified with partial lattices $P$ such that
$\bigvee$ and $\bigwedge$ are defined for all finite subsets of $P$. We
denote by $\Conc P$ the \jzs\ of all compact congruences of $P$. The
compact congruences of $P$ are those of the form
  \[
  \bigvee_{i<n}\Theta^+(a_i,b_i),
  \]
where $n<\omega$, $a_0$, \dots, $a_{n-1}$, $b_0$, \dots,
$b_{n-1}\in P$, and where we define $\Theta^+(a_i,b_i)$ to be the least
congruence $\preceq$ of $P$ such that $a_i\preceq b_i$. We observe
that congruences of a partial lattice $P$ are no longer equivalence
relations on $P$ but \emph{quasi-orderings} of~$P$.

\begin{definition}\label{D:HomPL}
If $P$ and $Q$ are partial prelattices, a \emph{homomorphism of partial
prelattices} from $P$ to $Q$ is an order-preserving map
$f\colon P\to Q$ such that $a=\bigvee X$
(resp., $a=\bigwedge X$) implies that $f(a)=\bigvee f[X]$ (resp.,
$f(a)=\bigwedge f[X]$), for all $a\in P$ and all $X\in\fine P$. We say
that a homomorphism $f$ is an \emph{embedding}, if $f(a)\leq f(b)$ implies
that $a\leq b$, for all $a$, $b\in P$.
\end{definition}

For a homomorphism $f\colon P\to Q$ of partial lattices the map
$\Conc f\colon\Conc P\to\Conc Q$ assigns to every compact congruence
$\xa$ of $P$ the congruence of $Q$ generated by all the pairs
$\seq{f(x),f(y)}$, for $\seq{x,y}\in\xa$. This way, the correspondence $\Conc$
becomes a functor from the category $\PL$ of partial lattices and
homomorphisms of partial lattices to the category $\CS$.

\section{Lifting diagrams of semilattices by diagrams of partial
lattices}\label{S:LiftPL}

In this section, we shall review some useful categorical concepts.

Every poset $I$ can be viewed as a category, whose set of objects is
$I$, where for all $p$, $q\in I$, there exists a morphism from $p$ to
$q$ exactly when $p\leq q$, and then this morphism is unique; we shall
denote it by $p\to q$. We shall often identify a poset with its
associated category.

For a category $\CC$, we shall say that a \emph{diagram} of $\CC$ is a
functor $\DD\colon I\to\CC$, where $I$ is a poset. If $\CD$ is another
category, a functor $\Phi\colon\CC\to\CD$ is said to \emph{preserve
direct limits} ( = directed colimits), if whenever $\DD$ is a diagram of
$\CC$ indexed by a \emph{directed} poset $I$ and $X=\varinjlim\DD$ in
$\CC$, then $\Phi(X)=\varinjlim(\Phi\DD)$. Observe that this
standard category-theoretical formulation abuses notation in two ways:
\begin{itemize}
\item Strictly speaking, $X$ does not simply consist of an object of
$\CC$ but rather of an object of $\CC$ together with a family of
morphisms $\DD(i)\to X$, for $i\in I$, satisfying natural commutation 
relations.

\item The statement $X=\varinjlim\DD$ determines $X$ only up to
isomorphism.
\end{itemize}

\begin{proposition}\label{P:PresLim}
The functor $\Conc$ from partial lattices to \jzs s preserves direct
limits.
\end{proposition}

Of course, Proposition~\ref{P:PresLim} is not specific to partial
lattices, it is an easy basic fact of universal algebra that holds for any
``reasonable'' definition of a partial algebra.

Now let us see how this can help us tackle \CLP.

Since every distributive semilattice $D$ is the direct union of
its finite distributive subsemilattices (see P. Pudl\'ak \cite{Pudl}),
one can start with a diagram
$\EE\colon I\to\CSfd$ indexed by a directed poset $I$ and try to find
a diagram $\DD\colon I\to\CL$ such that $\varinjlim\Conc\DD$ and
$\varinjlim\EE$ are isomorphic. This is true if the functors
$\Conc\DD$ and $\EE\colon I\to\CSfd$ are naturally equivalent (see
P. Pudl\'ak \cite{Pudl}), that is, if there exists a system
$\famm{\eps_i}{i\in I}$ of isomorphisms $\eps_i\colon\Conc\DD(i)\to\EE(i)$ (for
$i\in I$) such that the diagram of Figure~2 commutes, for all $i\leq j$ in $I$.
  \begin{figure}[hbt]
  \includegraphics{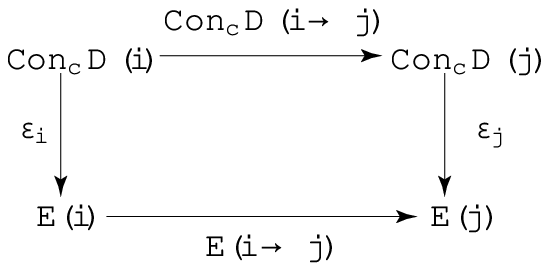}
  \caption{Natural equivalence of $\Conc\DD$ and $\EE$.}
  \end{figure}
If this is the case, then we also say that $\DD$ \emph{lifts} $\EE$ \emph{with
respect to} $\Conc$.

In some cases, lifts of diagrams can be constructed inductively using
the following concepts, see F. Wehrung
\cite{Wehr2con}. By a morphism $\varphi\colon\EE\to S$ of a diagram
$\EE\colon I\to\CS$ to an object $S$ of $\CS$ we mean a family
$\famm{\varphi_i}{i\in I}$ of morphisms
$\varphi_i\colon\EE(i)\to S$ of $\CS$ such that the equality
$\varphi_i=\varphi_j\circ\EE(i\to j)$ holds, for all $i\leq j$ in
$I$.

\begin{definition}\label{D:FactLift}
Let $\DD$ be a diagram of partial lattices.
For a \jzs\ $S$ and a partial lattice $P$, we say that a homomorphism
$\varphi\colon\Conc\DD\to S$ can be
\begin{enumerate}
\item \emph{factored through} $P$, if there are a homomorphism
$f\colon\DD\to P$ and a \jzh\ $\psi\colon\Conc P\to S$ such that
$\varphi=\psi\circ\Conc f$;

\item\emph{lifted through} $P$, if there are a homomorphism
$f\colon\DD\to P$ and an isomorphism $\psi\colon\Conc P\to S$ such that
$\varphi=\psi\circ\Conc f$.
\end{enumerate}

In (i) (resp., (ii)) above, we say that $\varphi$ can be \emph{factored
to} (resp., \emph{lifted to}) $f$.
\end{definition}

For example, any distributive semilattice $C$ is isomorphic to the
direct limit of a diagram $\EE\colon\fin{C\times\omega}\to\CSfb$ of
finite Boolean semilattices, see K.\,R. Goodearl and F.
Wehrung \cite{GoWe} ($\fin{X}$ denotes the set of all finite subsets
of a set $X$). So if this diagram could be lifted with respect to
the functor $\Conc$, then it would give a proof that every distributive
semilattice is isomorphic to the semilattice of compact congruences of
a lattice.

By a \emph{truncated $n$-cube} of lattices, we mean a functor
$\EE\colon\PP_<(n)\to\CL$, where we put
$\PP_<(n)=\PP(n)\setminus\set{n}$, partially ordered by inclusion. An
inductive construction of a lift (with respect to $\Conc$) of a
diagram $\EE\colon\fin{C\times\omega}\to\CSfb$ then requires
constructing lifts of homomorphisms of the form
$\varphi_n\colon\Conc\EE_n\to D$, where $n$ is a natural number, $D$ is
a distributive semilattice, and
$\EE_n\colon\PP_<(n)\to\CL$ is a truncated $n$-cube of lattices,
through a lattice $L$. This suggests that the study of the following
\emph{$n$-dimensional versions of \CLP} may be of interest.

\begin{definition}\label{D:nCLP}
For a natural number $n$, we say that a \jzs\ $S$ satisfies
\emph{$n$-dimensional \CLP}, or \emph{$n$-\CLP} in short, if for every
truncated $n$-cube $\EE$ of lattices, every homomorphism
$\varphi\colon\Conc\EE\to S$ can be lifted.
\end{definition}

Thus, $(n+1)$-\CLP\ at $S$ is stronger than $n$-\CLP\ at $S$, while 
$0$-\CLP\ at
$S$ is equivalent to the statement that $S\cong\Conc L$, for some 
lattice $L$. It
is worthwhile to restate $n$-\CLP\ at $S$, for $n\in\set{1,2}$:

\begin{itemize}
\item[\textbf{$1$-\CLP.}] The property $1$-\CLP\ holds at $S$ if{f} for every
lattice $K$, every \jzh\ $\varphi\colon\Conc K\to S$ can be lifted,
that is, there are a lattice $L$, a lattice homomorphism
$f\colon K\to L$, and an isomorphism $\eps\colon\Conc L\to S$ such
that $\varphi=\eps\circ\Conc f$, as illustrated on Figure~3.
  \begin{figure}[hbt]
  \includegraphics{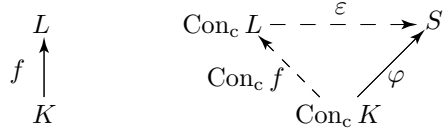}
  \caption{Illustrating $1$-\CLP\ at $S$.}
  \end{figure}
%
\item[\textbf{$2$-\CLP.}]
The property $2$-\CLP\ holds at $S$ if{f} for all lattices $K_0$,
$K_1$, $K_2$, all lattice homomorphisms $f_i\colon K_0\to K_i$ and all
\jzh s $\varphi_i\colon\Conc K_i\to S$, for $i\in\set{1,2}$ such that
$\varphi_1\circ\Conc f_1=\varphi_2\circ\Conc f_2$, there are a lattice~$L$,
lattice homomorphisms $g_i\colon K_i\to L$, for $i\in\set{1,2}$, and 
an isomorphism
$\eps\colon\Conc L\to S$ such that $g_1\circ f_1=g_2\circ f_2$ and
$\eps\circ\Conc g_i=\varphi_i$, for $i\in\set{1,2}$, as illustrated 
on Figure~4.
  \begin{figure}[hbt]
  \includegraphics{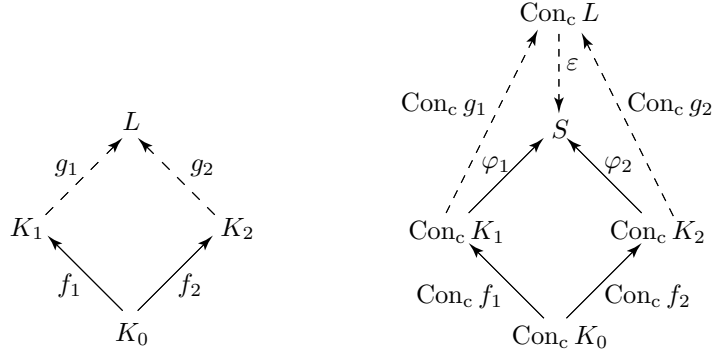}
  \caption{Illustrating $2$-\CLP\ at $S$.}
  \end{figure}
%
\end{itemize}

For $n\geq3$, it follows from a simple example in F. Wehrung
\cite{Wehr2con} that $n$-\CLP\ has a trivial answer:

\begin{proposition}\label{P:3CLP}
Let $S$ be a \jzs. Then $3$-\CLP\ holds at $S$ if{f} $S=\set{0}$.
\end{proposition}

A much harder related result is proved in J. T\r{u}ma and F. Wehrung
\cite{TuWe1}:

\begin{theorem}\label{T:CubeNL}
There exists a $3$-cube $\EE\colon\PP(3)\to\CSfb$ that has no lift
$\DD\colon\PP(3)\to\CL$ with respect to $\Conc$ such that
$\DD(\set{i})$ has \pup{almost} permutable congruences, for all $i<3$.
\end{theorem}

Here, we say that a lattice $L$ has \emph{permutable
congruences} (resp., \emph{almost permutable congruences}),
if $\xa\xb=\xb\xa$ (resp., $\xa\vee\xb=\xa\xb\cup\xb\xa$), for all
congruences $\xa$ and $\xb$ of $L$. The basic idea underlying the proof of
Theorem~\ref{T:CubeNL} consists of extracting the combinatorial core of
Theorem~\ref{T:PermConURPm}.

\section{Extensions of partial lattices to lattices}
\label{S:PartLat}

The main result of F. Wehrung \cite{Wehr1con}, together with the
converse proved in J.~T\r{u}ma and F. Wehrung \cite{TuWe2}, imply
the following.

\begin{theorem}\label{T:Con1}
Let $S$ be a distributive \jzs\ $S$. Then $1$-\CLP\ holds at~$S$ if{f}
$S$ is a lattice.
\end{theorem}

The proof that $1$-\CLP\ at $S$ implies that $S$ is a lattice is
established in J.~T\r{u}ma and F. Wehrung \cite{TuWe2}, and it uses
an \emph{ad hoc} construction: for a distributive \jzs\ $S$ that is
not a lattice, one constructs a \emph{Boolean} lattice $B$ of
cardinality $2^{|S|}$ and a \jzh\ $\varphi\colon\Conc B\to S$ without a lift.
\smallskip

The proof that $S$ being a lattice implies $1$-\CLP\ at $S$ is
very long and difficult, although the basic idea is quite simple. We are
given a distributive \jzs\ $S$, a lattice $K$, and a \jzh\
$\varphi\colon\Conc K\to S$. This information can be conveniently
expressed by saying that the pair $\seq{K,\varphi}$ is a \emph{\ML{S}}. We
extend $\seq{K,\varphi}$ to a pair $\seq{L,\psi}$, by successive one-step
extensions starting on $\seq{K,\varphi}$ that can be informally described
as follows:

\begin{enumerate}
\item Let $a<b<c$ in $K$. We freely adjoin to $K$ a relative
complement $x$ of $b$ in the interval $[a,c]$. Then the map $\psi$
sends $\Theta_L(a,x)$ to $\varphi\Theta_K(b,c)$ and $\Theta_L(x,c)$ to
$\varphi\Theta_K(a,b)$.

\item Let $\xa$ and $\xb$ be congruences of $K$ such that
$\varphi(\xa)=\varphi(\xb)$. Then $L$ is constructed such
that $(\Conc j)(\xa)=(\Conc j)(\xb)$, where $j\colon K\to L$ is a
lattice homomorphism, while the construction is ``sufficiently free''.

\item Let $\alpha\in S$ be not in the range of $\varphi$, fix $o\in K$
and add freely to $K$ an outside element $x>o$. Let $\psi$ send
$\Theta_L(o,x)$ to $\alpha$.
\end{enumerate}

Iterating these steps transfinitely should, intuitively, yield an \ML{S}
$\seq{L,\psi}$ with $L$ relatively complemented and $\psi$ an 
\emph{isomorphism}.
This approach suffers from many drawbacks:

\begin{itemize}
\item[(a)] The structure $L$ obtained above is not a lattice, but a 
\emph{partial
lattice}. Hence the induction step should be performed not on a pair
$\seq{K,\varphi}$ with $K$ a \emph{lattice}, but on a pair $\seq{P,\varphi}$,
where $P$ is a \emph{partial lattice} and $\varphi\colon\Conc P\to S$
is a \jzh. We shall say that $\seq{P,\varphi}$ is a \emph{\MPL{S}}.

\item[(b)] Performing the steps above on an \MPL{S} $\seq{P,\varphi}$ leads
to the problem of extending the map $\varphi$ that may not have a
solution even in simple cases.

\item[(c)] Once the inductive construction is completed, taking the
limit does not yield a lattice, but a partial lattice.
\end{itemize}

Item (a) is taken care of by
the extension of the $\Conc$ functor to the category~$\PL$ of partial
lattices, as presented in Section~\ref{S:Conc}.

Item (b) is much more difficult to take care of. In F.
Wehrung \cite{Wehr1con}, this is done by considering a certain class
of \MPL{S}s $\seq{P,\varphi}$ that are called \emph{balanced} \MPL{S}s, for
which the extension from $\varphi$ to $\psi$ can be performed. The
formal definition of being balanced is quite complicated. Intuitively,
it means that the meet and join operations can be computed in the ideal
and filter lattices of $P$ by focusing attention on
\emph{finite} sets of elements, called there \emph{samples}, and do
this uniformly on all quotients of $P$.

Item (c) is taken care in a similar fashion as item (b). For a partial lattice~$P$,
let $\FL(P)$ denote the free lattice on $P$, see R.\,P. Dilworth 
\cite{Dilw45},
R.\,A. Dean \cite{Dean64}, or R. Freese, J. Je\v{z}ek, and J.\,B. 
Nation \cite{FJN}.
Then the canonical map from $\Conc P$ to $\Conc\FL(P)$ is a cofinal
\jz-embedding. This alone is not sufficient to ensure the existence of
an extension $\psi\colon\Conc\FL(P)\to S$ of the map
$\varphi\colon\Conc P\to S$, however, this is possible if $\seq{P,\varphi}$
is balanced. The formula defining the extension $\psi$ may be better
understood by viewing $\seq{P,\varphi}$ as a \emph{$S^\dd$-valued partial
lattice}, where~$S^\dd$ denotes the dual lattice of $S$ (thus it is a
distributive lattice with $1$). Thus we need to deal with structures
similar to the Boolean-valued models encountered in forcing
(e.g., in set theory), except that the ``truth values'' live
not in a complete Boolean algebra but in the lattice $S^\dd$.

Once all these problems are solved, we obtain the following
much stronger result, see Theorem~D in F. Wehrung \cite{Wehr1con}:

\begin{theorem}\label{T:1DimLift}
Let $S$ be a distributive lattice with zero. Then for all lattices
$K_0$, $K_1$, $K_2$ with $K_0$ \emph{finite}, all lattice homomorphisms
$f_i\colon K_0\to K_i$ and all \jzh s $\varphi_i\colon\Conc K_i\to S$,
for $i\in\set{1,2}$ such that
$\varphi_1\circ\Conc f_1=\varphi_2\circ\Conc f_2$, there are a lattice
$L$, lattice homomorphisms $g_i\colon K_i\to L$, for $i\in\set{1,2}$,
and an isomorphism $\eps\colon\Conc L\to S$ such that
$\eps\circ\Conc g_i=\varphi_i$, for $i\in\set{1,2}$. Furthermore, $L$,
$g_1$, and $g_2$ can be found in such a way that the following
statements hold:

\begin{enumerate}
\item $L$ is relatively complemented.

\item $g_1[K_1]\cup g_2[K_2]$ generates $L$ as an ideal \pup{resp., a
filter}.

\item If $\rng\varphi_1\cup\rng\varphi_2$ generates $S$ as an ideal,
then $g_1[K_1]\cup g_2[K_2]$ generates $L$ as a convex sublattice.
\end{enumerate}
\end{theorem}

This is a far reaching generalization of the original result of M. Tischendorf
\cite{Tisch} with $K_0=K_1=K_2$ finite, $f_1=f_2=\id$, $S$ finite, and
$\varphi$ an embedding.

The lifts satisfying conditions (i)--(iii) in 
Theorem~\ref{T:1DimLift} are called
\emph{good lifts} in F.~Wehrung \cite{Wehr1con}. The properties of
$\seq{L,\eps}$ from which they follow are reminiscent of \emph{genericity}
(in the model-theoretical sense), and these properties have further
consequences, for example, the lattice $L$ has \emph{definable
congruence inclusion} in the sense that for every positive integer $n$, there
exists a positive existential formula (independent of the lattice $L$)
$\Phi_n(\mathsf{x}_0,\mathsf{y}_0,\mathsf{x}_1,\mathsf{y}_1,
\ldots,\mathsf{x}_n,\mathsf{y}_n)$
of lattice theory such that $L$ satisfies that
  \[
  \Theta(x_0,y_0)\subseteq\bigvee_{i=1}^n\Theta(x_i,y_i)\qquad\text{if{f}}\qquad
  L\text{ satisfies }\Phi_n(x_0,y_0,x_1,y_1,\ldots,x_n,y_n),
  \]
for any $x_0$, $y_0$, \dots, $x_n$, $y_n\in L$.
As for the lifts of truncated squares (the case $n=2$) earlier results
by the first author, J. T\r{u}ma \cite{Tuma} and also G.~Gr\"atzer,
H.~Lakser, and F.~Wehrung~\cite{GLW} are extended to infinite
semilattices in F.~Wehrung \cite{Wehr2con}.

\begin{definition}\label{D:CCB}
A \jzs\ $S$ is called \emph{conditionally co-Brouwerian}, if the 
following holds:
\begin{enumerate}
\item for all nonempty subsets $X$ and $Y$ of $S$ such that $X\leq Y$
(that is, $x\leq y$, for all $x\in X$ and $y\in Y$), there exists $z\in S$ such
that $X\leq z\leq Y$ (we then say that $S$ is \emph{conditionally complete});

\item for every subset $Z$ of $S$, if $a\leq b\vee z$, for all
$z\in Z$, then there exists $c\in S$ such that $a\leq b\vee c$ and
$c\leq Z$.
\end{enumerate}

By restricting this definition to subsets $X$, $Y$, and $Z$
of $S$ such that $|X|$, $|Y|$, $|Z|<\kappa$, for an infinite cardinal
$\kappa$, we define \emph{\ckcb} \jzs s.
\end{definition}

Of course, any \ccb\ \jzs\ is a distributive lattice with zero.

\begin{theorem}\label{T:2Main}
Let $\EE$ be a truncated square of partial lattices with
$\EE(\es)$ a lattice, let~$S$ be a \ccb\ lattice. Then every
homomorphism $\varphi\colon\Conc\EE\to S$ has a good lift.
\end{theorem}

A self-contained proof of Theorem~\ref{T:2Main} is significantly
easier than a self-contained proof of Theorem~\ref{T:1DimLift},
because for a \ccb\ lattice $S$ and a partial lattice $P$, any \jzh\
$\varphi\colon\Conc P\to S$ can be extended to a \jzh\
$\psi\colon\Conc\FL(P)\to S$; this follows from
monoid-theoretical considerations introduced in F.
Wehrung~\cite{Wehr92}. The assumption that $S$ is \ccb\ implies that
$S$ is \emph{injective} in a suitable category of partially
quasi-ordered monoids.

As a consequence of one- and two-dimensional lifting results in F.
Wehrung \cite{Wehr1con,Wehr2con}, we obtain the following
extensions of known $0$-dimensional results:

\begin{theorem}\label{T:New0Dim}
Every member in each of the following classes of \jzs s is isomorphic to
$\Conc L$, for some relatively complemented lattice $L$ with zero that
has definable congruence inclusion:
\begin{enumerate}
\item \jz-direct limits of the form $\varinjlim_{n\in\omega}S_n$, where
all the $S_n$ are distributive \emph{lattices} with zero;

\item \jz-direct limits of the form $\varinjlim_{i\in I}S_i$, where
$I$ is a directed poset of cardinality at most $\aleph_1$ and all the
$S_i$, for $i\in I$, are \ccb;

\item All \jzs s $S$ that are conditionally $|S|$-co-Brou\-wer\-ian.
\end{enumerate}
\end{theorem}

Item (i) above extends the main result of E.\,T. Schmidt \cite{Schm81}
(any distributive lattice with zero is representable), while (ii)
above extends the main result of A.\,P. Huhn \cite{Huhn89a,Huhn89b} (any
distributive \jzs\ of cardinality at most $\aleph_1$ is representable).
Item (iii), first stated and proved in F. Wehrung \cite{Wehr2con}, seems to be
completely new.

We recall here the following result, obtained by E.\,T. Schmidt, see 
\cite{Schm68}:

\begin{theorem}[Schmidt's Lemma]\label{T:Schmidt}
Let $B$ be a generalized Boolean semilattice. Then any image of $B$ under a
distributive homomorphism \pup{see Definition~\textup{\ref{D:DistrCon}}} is
isomorphic to $\Conc L$, for some lattice $L$.
\end{theorem}

We obtain the following informal corollary:

\begin{corollary}\label{C:Informal}
All the representation theorems of distributive \jzs s other than 
Schmidt's Lemma
that are known to this point yield semilattices that are representable by
relatively complemented lattices with zero and with definable 
congruence inclusion.
\end{corollary}

It is unclear whether every image $S$ of a generalized Boolean 
semilattice under a
distributive homomorphism can be represented as $\Conc L$, for a relatively
complemented lattice $L$. However, a direct verification yields that 
$S$ satisfies
the axiom $\URP$ considered in Section~\ref{S:URP} (see
Proposition~\ref{P:BasicWD}(ii)), and the (easy) proof fails for the 
stronger axiom
$\URP^+$ (see Definition~\ref{D:URPp}). On the other hand, the 
semilattice $\Conc L$
satisfies $\URP^+$, for any relatively complemented lattice $L$ (see
Proposition~\ref{P:RelCplURP+}); see also Problem~\ref{Pb:URPSchm}.

In any case, all known representation theorems (Schmidt's Lemma included) yield
semilattices that satisfy the axiom $\URP$ studied in Section~\ref{S:URP}. In
particular, none of them is able to reach $\Conc\FL(\omega_2)$, although this
semilattice is already represented!

\section{Connections to ring theory}\label{S:Ring1}

The paper K.\,R. Goodearl and F. Wehrung \cite{GoWe} is a rich source of
information on connections between congruence lattice representation
problems and ring theory. For our present purpose, we mention the
following theorem that goes back to J.~von~Neumann \cite{vN}. A ring
(associative, not necessarily with unit) $R$ is called \emph{regular}
(in von~Neumann's sense), if for all
$x\in R$, there exists $y\in R$ such that $xyx=x$. We recall the
following classical result, see K.\,D. Fryer and I. 
Halperin~\cite{FrHa56} for a
proof of the case without unit:

\begin{proposition}\label{P:rrings}
The set $\LL(R)$ of all principal right ideals of a regular ring
$R$, ordered under inclusion, is a sectionally complemented modular
lattice.
\end{proposition}

This together with the following result from F. Wehrung \cite{WeURP}
gives a strategy for representing distributive semilattices as
semilattices of compact congruences of sectionally complemented modular
lattices.

\begin{proposition}\label{P:rrings2}
Let $R$ be a regular ring. Then the semilattices $\Conc\LL(R)$ and
$\Idc R$ \pup{the semilattice of finitely generated two-sided ideals of
$R$} are isomorphic \pup{distributive} semilattices.
\end{proposition}

\emph{A matricial algebra} over a field $F$ is a finite direct product
of full matricial algebras over $F$. A direct limit of matricial
algebras over
$F$ is called {a locally matricial algebra} over $F$. Locally matricial
algebras are regular. If $F$ is a \emph{finite} field and $R$ is a
locally matricial algebra over $F$, then $R$ is a locally finite ring
and $\LL(R)$ is a locally finite lattice.
The following theorem appears first in well-known unpublished
notes of G.\,M. Bergman \cite{Berg86}:

\begin{theorem}\label{T:berg}
Let $F$ be a field. Then every countable distributive \jzs\ is isomorphic to
$\Idc R$, for some locally matricial algebra $R$ over $F$.
\end{theorem}

P. R\r{u}\v{z}i\v{c}ka proves in \cite{Ruzi} by a very sophisticated
construction that Bergman's theorem holds also for distributive
\emph{lattices} of arbitrary cardinality.

\begin{theorem}\label{T:Ruz}
Let $F$ be a field. Then every distributive lattice with zero is
isomorphic to $\Idc R$, for some locally matricial algebra $R$ over $F$.
\end{theorem}

The problem whether every distributive semilattice of cardinality
$\aleph_1$ is isomorphic to $\Idc R$, for some locally matricial ring
$R$, remains open, see Problem~\ref{Pb:Al1Nab}. Because of the
following result in F. Wehrung \cite{We3}, the cardinality
$\aleph_1$ is the maximal cardinality for which there could be a general
positive answer:

\begin{theorem}\label{T:Norrings}
There exists a distributive semilattice of cardinality $\aleph_2$ that
is not isomorphic to $\Idc R$, for any regular ring $R$.
\end{theorem}

On the other hand, every distributive \jzs\ of cardinality
$\aleph_1$ can be represented as $\Idc R$, for some regular ring $R$, as
also proved in F. Wehrung \cite{We4}. In the spirit of
Problem~\ref{Pb:Al1Nab}, we can also mention the following difficult
one-dimensional analogue of Theorem~\ref{T:berg}, established in
J. T\r{u}ma and F. Wehrung \cite{TuWe3}:

\begin{theorem}\label{T:LiftCtbleS2T}
Let $F$ be a field, let $S$ and $T$ be countable distributive \jzs s,
let $\varphi\colon S\to T$ be a \jzh. Then there are
locally matricial algebras $A$ and $B$ over $F$, a
homomorphism $f\colon A\to B$ of $F$-algebras, and isomorphisms
$\alpha\colon\Idc A\to S$ and $\beta\colon\Idc B\to T$ such that
$\beta\circ\Idc f=\varphi\circ\alpha$.
\end{theorem}

An interesting point about the proof of Theorem~\ref{T:LiftCtbleS2T} is that it
involves a \emph{reverse one-dimensional amalgamation} result. We say that a
partially ordered vector space (over the field $\QQ$ of rational numbers) is
\emph{simplicial}, if it is isomorphic to a finite power of $\QQ$ 
with componentwise
ordering, and that it is a \emph{dimension vector space}, if it is 
isomorphic to a
direct limit of simplicial vector spaces. The main result in J.~T\r{u}ma and F.
Wehrung~\cite{TuWe3} is that for a countable dimension vector space $V$ and a
countable distributive \jzs~$S$, every \jzh\ from $S$ to $\Idc V$ can 
be lifted by a
positive homomorphism from $U$ to $V$, for some
(countable) dimension vector space $U$. The final step from dimension groups to
locally matricial algebras uses results from K.\,R. Goodearl and D.\,E.
Handelman~\cite{GoHa}.

Theorem~\ref{T:LiftCtbleS2T} yields the following lattice-theoretical 
consequence:

\begin{corollary}\label{C:LiftCtbleS2T}
Let $S$ and $T$ be countable distributive \jzs s, let\linebreak
$\varphi\colon S\to T$ be a \jzh. Then there are locally finite,
relatively complemented modular lattices $K$ and $L$, a
lattice homomorphism $f\colon K\to L$, and isomorphisms
$\alpha\colon\Conc K\to S$ and $\beta\colon\Conc L\to T$ such that
$\beta\circ\Conc f=\varphi\circ\alpha$.
\end{corollary}

It is also established in J. T\r{u}ma and F. Wehrung \cite{TuWe3} that
the analogue of Corollary~\ref{C:LiftCtbleS2T} for $S$
\emph{uncountable} fails.

\section{Dual topological spaces}\label{S:AlgLatt}

In the important papers \cite{Plo1} and \cite{Plo2},
M.~Plo\v{s}\v{c}ica investigates dual topological spaces of some
congruence lattices. Any algebraic distributive lattice $D$ defines a
topological space $\M(D)$. The points of $\M(D)$ are completely
meet-irreducible elements of $D$ and closed sets of $\M(D)$ are sets of
the form $\M(D)\cap[x,1_D]$, for any $x\in X$. The lattice $D$ can be
reconstructed from its dual space
$\M(D)$ as the lattice of open subsets of $\M(D)$ ordered by inclusion.

If $L$ is a lattice, then the points of the dual space $\M(\Con L)$ of
the full congruence lattice of $L$ are the subdirectly irreducible
congruences of $L$, that is, the congruences $\xa$ of $L$ such
that the quotient lattice $L/{\xa}$ is subdirectly irreducible. It
seems that the dual spaces $\M(\Con L)$ might be a useful tool in the
study of the congruence lattices of members of lattice varieties with
only finitely many non-isomorphic subdirectly irreducible lattices.

The dual spaces $\M(D)$ have a base of compact open sets but they are
not usually Hausdorff.

\begin{definition}\label{D:ConCl}
For a class $\mathbf{C}$ of lattices, we define
  \[
  \Con\mathbf{C}=\setm{D}{D\cong\Con L,\text{ for some }L\in\mathbf{C}},
  \]
the \emph{congruence class} of $\mathbf{C}$.
\end{definition}

M. Plo\v{s}\v{c}ica proves in \cite{Plo1} that the
congruence classes $\Con\MM_n$ (here $\MM_n$ denotes the variety
generated by $M_n$, the lattice of length two with $n+2$ elements) are
distinct. The topological property that distinguishes them is
\emph{uniform separability}.

\begin{definition}\label{D:Separ}
A subset $Q$ of a topological space $T$ is called \emph{discrete}, if
every subset of $Q$ is open in the relative topology on $Q$.
The space $T$ is called \emph{uniformly}
$n$\emph{-separable} (for $n\geq 3$), if for every discrete set $Q\subseteq T$,
there exists a family $\famm{U_{pq}}{p,\,q\in Q,\ p\ne q}$ of open 
sets such that
$p\in U_{pq}$, for every $p$, $q\in Q$, and, for every $n$-element set
$Q_0\subseteq Q$,
  \[
  \bigcap\setm{U_{pq}}{p,q\in Q_0,\ p\ne q}=\es.
  \]
\end{definition}

The following two theorems establish the crucial separability
properties of the spaces $\M(\Con L)$, for $L\in\MM_n$.

\begin{theorem}\label{T:Separ1}
If $L\in\MM_n$, $n\geq 3$, then $\M(\Con L)$ is $(n+1)$-uniformly
separable.
\end{theorem}

In order to prove Theorem~\ref{T:Separ1}, M. Plo\v{s}\v{c}ica 
\cite{Plo1} assumes
that $(n+1)$-uniform separability fails in $\M(\Con L)$, and infers, with the
help of a clever combinatorial statement, that $M_n$ has $n+1$ 
distinct atoms, a
contradiction.

Let $F_n(X)$ denote the free lattice over $X$ in the variety $\MM_n$, 
for $n\geq 3$
and any set $X$.

\begin{theorem}\label{T:Separ2}
The topological space $\M(\Con F_n(X))$ is not $n$-uniformly separable.
\end{theorem}

As the proof of Theorem~\ref{T:PermConURP}(ii) is based on the 
Kuratowski Free Set
Theorem, the proof of Theorem~\ref{T:Separ2} is based on the 
following extension of
that theorem, due to A. Hajnal and A. M\'at\'e \cite{HaMa75}:

\begin{theorem}\label{T:HaMa}
Let $X$ be a set of cardinality at least $\aleph_2$, let
$\Phi\colon[X]^2\to\fin{X}$. Then for every natural number $n\geq3$, 
there exists
$U\in[X]^n$ such that $u\notin\Phi(V)$, for all $u\in U$
and all $V\in[U\setminus\set{u}]^2$.
\end{theorem}

As a corollary, we get the following:

\begin{corollary}\label{C:Separ3}
Let $n\geq 3$, let $X$ be a set of cardinality at least $\aleph_2$. Then there
is no lattice $L\in\MM_n$ such that $\Con L$ is isomorphic to $\Con 
F_{n+1}(X)$.
\end{corollary}

In his other paper \cite{Plo2}, M. Plo\v{s}\v{c}ica characterizes dual
spaces $\M(\Con L)$, for lattices~$L$ with at most $\aleph_1$ compact
elements from the variety $\MM_n^{01}$ generated by~$M_n$ as a bounded
lattice, $n\geq 3$. His main result is the following deep theorem:

\begin{theorem}\label{T:Separ4}
Let $D$ be an algebraic distributive lattice with at most $\aleph_1$
compact elements and $n\geq 3$. Then $D$ is isomorphic to $\Con L$, for some
$L\in\MM_n^{01}$, if and only if the topological space $T=\M(D)$ has a
subspace~$T_0$ that satisfies the following five conditions:
  \begin{enumerate}
  \item $T$ is compact and has a basis of compact open sets;
  \item both $T_0$ and $T_n=T\setminus T_0$ are Hausdorff zero-dimensional;
  \item $T_0$ is a closed subspace of $T$;
  \item if $a\in T_n$, $b\in T\setminus\set{a}$, then there exists a
  clopen set $V\subseteq T_n$ such that $a\in V$ and $b\notin V$;
  \item if $a$, $b$, $c\in T$ are distinct, then there exist
  open sets $U$, $V$, $W$ such that $a\in U$, $b\in V$, $c\in W$, and
  $U\cap V\cap W=\es$.
  \end{enumerate}
\end{theorem}

In order to establish the harder direction of Theorem~\ref{T:Separ4}, M.
Plo\v{s}\v{c}ica embeds directly, \emph{via} an elaborate \emph{ad hoc}
construction, the space $T$ as a closed subspace of 
$\M(\Con(F_n^{01}(\omega_1)))$,
where $F_n^{01}(\omega_1)$ denotes the free object on $\aleph_1$ generators in
the variety $\MM_n^{01}$.

Since the conditions on $\M(D)$ do not depend on $n\geq 3$, we get the
following corollary, see M. Plo\v{s}\v{c}ica \cite{Plo2}.

\begin{corollary}\label{C:Separ5}
If $L\in\MM_n^{01}$, $n\geq 3$, and $L$ has at most $\aleph_1$ 
elements, then there
exists $K\in\MM_3^{01}$ such that $\Con K$ is isomorphic to $\Con L$.
\end{corollary}

These results further emphasize the crucial role that the cardinality
$\aleph_2$ plays in the study of congruence lattices of lattices; see also
Problem~\ref{Pb:CritPt}.

\section{Open problems}\label{S:Pbs}

There are many open problems related to \CLP, scattered in the
literature. Here are, to our minds,
the most outstanding ones.

We first restate Dilworth's still unsolved problem:

\begin{all}{Congruence Lattice Problem}
Let $S$ be a distributive \jzs. Does there exist a
lattice $L$ such that $\Conc L\cong S$?
\end{all}

Our next problem is a byproduct of the study of $1$-\CLP:

\begin{problem}\label{Pb:CtbleCon1}
Let $K$ be a countable lattice, let $S$ be a countable distributive
\jzs, let $\varphi\colon\Conc K\to S$ be a \jzh. Can $\varphi$ be
lifted, that is, are there a lattice $L$, a lattice homomorphism
$f\colon K\to L$, and an isomorphism $\eps\colon\Conc L\to S$ such
that $\eps\circ\Conc f=\varphi$?
\end{problem}

We observe that the cardinality assumption on $K$ and $S$ in
Problem~\ref{Pb:CtbleCon1} is \emph{optimal}. Indeed, the paper
J. T\r{u}ma and F. Wehrung \cite{TuWe2} contains an example of a \jzh\
$\varphi\colon\Conc B\to S$, where $B$ is a Boolean lattice of
cardinality $\aleph_1$ and $S$ is a countable distributive \jzs, that
cannot be lifted. We also observe that the analogue of
Problem~\ref{Pb:CtbleCon1}, where $K$ is only a \emph{partial lattice}, fails,
because of some results in F. Wehrung \cite{Wehr2con}. Nevertheless, we
still conjecture that Problem~\ref{Pb:CtbleCon1} has a positive
solution.

By Schmidt's Lemma (see Theorem~\ref{T:Schmidt}), every distributive image of a
generalized Boolean semilattice is representable. This suggests the following
problem:

\begin{problem}\label{Pb:WDIm}
Let $K$ be a lattice, let $S$ be a distributive \jzs, let
$\varphi\colon\Conc K\onto S$ be a surjective distributive \jzh. Can
$\varphi$ be lifted?
\end{problem}

Of course, if $\varphi$ can be lifted, then the semilattice $S$ is 
representable.
Problem~\ref{Pb:WDIm} is first stated in J. T\r{u}ma and F. Wehrung
\cite[Problem~2]{TuWe2}.

\begin{problem}\label{Pb:Al1Nab}
Is it the case that every distributive \jzs\ of cardinality~$\aleph_1$
is isomorphic to $\Conc L$, for some sectionally complemented, modular,
\emph{locally finite} lattice $L$?
\end{problem}

A ring-theoretical equivalent to Problem~\ref{Pb:Al1Nab} is whether
every distributive \jzs\ of cardinality
$\aleph_1$ is isomorphic to $\Idc R$, for some locally matricial algebra
$R$, see K.\,R. Goodearl and F. Wehrung \cite{GoWe}. The countable case
is solved by Bergman's Theorem, see Theorem~\ref{T:berg}. The
statement obtained by removing ``locally finite'' from the statement
of Problem~\ref{Pb:Al1Nab} is proved in F. Wehrung~\cite{We4}. The
statement obtained by removing ``modular'' from the statement
of Problem~\ref{Pb:Al1Nab} is proved in G. Gr\"atzer, H. Lakser, and F.
Wehrung~\cite{GLW}. Various aspects and possible attacks of
Problem~\ref{Pb:Al1Nab} are also studied in J. T\r{u}ma and F. Wehrung
\cite{TuWe3}.

A problem related to Problem~\ref{Pb:CtbleCon1} is the following:

\begin{problem}\label{Pb:CPERC}
Does every lattice of cardinality at most $\aleph_0$ (resp., $\aleph_1$) have a
con\-gru\-ence-preserving extension to a relatively complemented lattice?
\end{problem}

It is proved in G. Gr\"atzer and E.\,T. Schmidt \cite{GrSc99} that every
finite lattice has a finite, sectionally complemented congruence-preserving
extension. Further results imply that every lattice $L$ in each of 
the following
classes has a relatively complemented congruence-preserving extension that it
generates as a convex sublattice:

\begin{itemize}
\item $L$ is a direct union $\bigcup_{n<\omega}L_n$, where $\Conc 
L_n$ is finite,
for all $n<\omega$ (G.~Gr\"atzer, H.~Lakser, and F.~Wehrung \cite{GLW}).

\item $L$ is a direct union $\bigcup_{n<\omega}L_n$, where $\Conc 
L_n$ is \ccb, for
all $n<\omega$ (F. Wehrung \cite{Wehr2con}).

\item $\Conc L$ is a lattice (F. Wehrung \cite{Wehr1con}).
\end{itemize}

On the other hand, the cardinality $\aleph_1$ in the statement of
Problem~\ref{Pb:CPERC} is the highest possible, because of the results
of M. Plo\v{s}\v{c}ica, J. T\r{u}ma, and F. Wehrung \cite{PTW} and
J. T\r{u}ma and F. Wehrung \cite{TuWe1}. For example, for any
nondistributive variety~$\VV$ of lattices, the free lattice in
$\VV$ on $\aleph_2$ generators does not have a
congruence-preserving extension with permutable (or even almost
permutable) congruences.

For varieties $\UU$ and $\VV$ of lattices, define
the \emph{critical point} of $\UU$ and $\VV$ as the least cardinality of the
semilattice of compact elements of a member of the symmetric difference
$(\Con\UU)\mathbin{\triangle}(\Con\VV)$ (let it be $\infty$, if 
$\Con\UU=\Con\VV$).

\begin{problem}[Critical point conjecture]\label{Pb:CritPt}
Let $\UU$ and $\VV$ be varieties of lattices (resp., finitely 
generated varieties
of lattices), with critical point $\kappa<\infty$. Prove that either
$\kappa\leq\aleph_0$ or $\kappa=\aleph_2$.
\end{problem}

In all known cases, the answer to Problem~\ref{Pb:CritPt} is positive. Even for
finitely generated varieties, this problems seems to be very difficult.

For a lattice $L$, let $\Var L$ denote the lattice variety generated by $L$,
and $\Con\Var L$ its congruence class (see Definition~\ref{D:ConCl}). 
We do not know
whether, for finite lattices $A$ and~$B$, the equality $\Con\Var 
A=\Con\Var B$ can
be checked recursively. In particular, the following sounds plausible:

\begin{problem}\label{L:DecVar}
For finite lattices $A$ and $B$, does $\Con\Var A=\Con\Var B$ imply that either
$A\cong B$ or $A\cong B^\dd$?
\end{problem}

It follows from M. Plo\v{s}\v{c}ica's results \cite{Plo1} that the
congruence classes $\Con\MM_n$, for $n\geq\nobreak 3$, are distinct, see
Section~\ref{S:AlgLatt}. However, they can be separated only by semilattices of
cardinality at least $\aleph_2$, see Corollary~\ref{C:Separ5}.

\begin{problem}\label{Pb:CharMn}
Characterize the congruence classes of the varieties $\MM_n$, for
$n\geq 3$.
\end{problem}

Up to now, the only nontrivial variety of which the congruence class is
completely described is the variety of \emph{distributive lattices},
for which the congruence class is the class of all lattices of ideals
of generalized Boolean semilattices (this is trivial and well-known). Necessary
conditions for a given distributive \jzs\ to belong to the congruence class of
$\MM_n$ are given in Theorems~\ref{T:Separ1} and \ref{T:Separ4}.

\begin{problem}\label{Pb:PCCC}
Are the congruence classes of sectionally complemented and relatively
complemented lattices distinct?
\end{problem}

Of course, further variants of Problem~\ref{Pb:PCCC} could be stated
for other classes of lattices, for example, the class of sectionally
complemented modular lattices or the class of lattices with 
permutable congruences.
A basic approach for tackling Problem~\ref{Pb:PCCC} is the following:

\begin{problem}\label{Pb:URPsc}
For a sectionally complemented lattice $L$, does $\Conc L$ satisfy
$\URP^+$?
\end{problem}

We proved in Proposition~\ref{P:RelCplURP+} that if $L$ is relatively
complemented, then $\Conc L$ satisfies $\URP^+$, but the argument
fails for sectionally complemented lattices.

\begin{problem}\label{Pb:URPSchm}
Let $S$ be a distributive \jzs\ that is the image of a generalized
Boolean semilattice under a distributive homomorphism. Does $S$ satisfy
$\URP^+$?
\end{problem}

\begin{problem}\label{Pb:URPTens}
Is the property $\URP$ preserved under tensor product (resp., iterated
tensor product) of \jzs s?
\end{problem}

Finally we observe that the problem about which semilattices $S$
satisfy $1$-\CLP\ (see Definition~\ref{D:nCLP}) is completely solved in
F. Wehrung \cite{Wehr1con} and J. T\r{u}ma and F.~Wehrung \cite{TuWe2}:
namely, these are exactly the distributive lattices with zero. The
two-dimensional analogue $2$-\CLP\ is not completely solved yet:

\begin{problem}\label{Pb:Con2}
Let $S$ be a distributive \jzs. Does any of the following assumptions
imply that $S$ is \ccb:
\begin{enumerate}
\item For every partial lattice $P$, every \jzh\
$\varphi\colon\Conc P\to S$ can be factored through a lattice;

\item For every partial lattice $P$, every \jzh\
$\varphi\colon\Conc P\to S$ can be lifted through a lattice;

\item For every truncated square $\DD$ of lattices, every \jzh\linebreak
$\varphi\colon\Conc\DD\to\nobreak S$ can be factored through a lattice;

\item For every truncated square $\DD$ of lattices, every \jzh\linebreak
$\varphi\colon\Conc\DD\to\nobreak S$ can be lifted through a lattice.
\end{enumerate}
\end{problem}

It is proved in F. Wehrung \cite{Wehr2con} that the assumption that $S$ be
\ccb\ is sufficient to imply (i)--(iv) above, see
Theorem~\ref{T:2Main}. Moreover, some partial converses to this
statement are proved in F. Wehrung
\cite{Wehr2con}: namely, either (i), (ii), or (iv) implies that $S$ is a
conditionally complete lattice.

\section*{Acknowledgments}

This work was partially completed while the second author was visiting
the Charles University (Prague). Excellent conditions provided by the
Department of Algebra are greatly appreciated. The authors also wish to thank
Marina Semenova and George Gr\"atzer for their comments.

\section*{Added in proof}
The second author recently solved Problem~\ref{Pb:Al1Nab}.

\end{document}